\documentclass[final]{siamltex}
\usepackage{amsmath}
\usepackage{amssymb}
\usepackage{url}
\usepackage{color}
\usepackage{array}
\usepackage{colortbl}
\usepackage{blkarray}
\usepackage{caption}
\usepackage{capt-of}
\usepackage{colortbl}
\usepackage{fancyhdr}
\usepackage{fancyvrb}
\usepackage{float}
\usepackage[T1]{fontenc}
\usepackage{framed}
\usepackage{geometry}
\usepackage{graphicx}
\usepackage[latin1]{inputenc}
\usepackage{stmaryrd}
\usepackage{upquote}
\usepackage[usenames,dvipsnames,svgnames,table]{xcolor}
\usepackage{xspace}
\usepackage{todonotes}
\usepackage{booktabs}
\definecolor{lightgray}{gray}{0.5}
\definecolor{darkgreen}{rgb}{0,0.6,0.13}
\newcommand{\nc}{\newcommand}
\nc{\dsp}{\displaystyle}
\nc{\txt}{\textstyle}
\nc{\reff}[1]{(\ref{#1})}
\nc{\mrm}[1]{\mathrm{#1}}
\nc{\udl}[1]{\underline{#1}}
\nc{\ovl}[1]{\overline{#1}}
\nc{\al}{\underline{\boldsymbol{\alpha}}}
\nc{\la}{\underline{\boldsymbol{\lambda}}}
\nc{\llbr}{\llbracket}
\nc{\rrbr}{\rrbracket}
\nc{\lbr}{\lbrack}
\nc{\rbr}{\rbrack}
\nc{\N}{\mathbb{N}}
\nc{\Z}{\mathbb{Z}}
\nc{\D}{\mathbb{D}}
\nc{\Q}{\mathbb{Q}}
\nc{\R}{\mathbb{R}}
\nc{\C}{\mathbb{C}}
\nc{\T}{\mathbb{T}}
\nc{\Abf}{\mathbf{A}}
\nc{\Bbf}{\mathbf{B}}
\nc{\Cbf}{\mathbf{C}}
\nc{\Ibf}{\mathbf{I}}
\nc{\Lbf}{\mathbf{L}}
\nc{\Nbf}{\mathbf{N}}
\nc{\Sbf}{\mathbf{S}}
\nc{\Mbf}{\mathbf{M}}
\nc{\Tbf}{\mathbf{T}}
\nc{\Dbf}{\mathbf{D}}
\nc{\Qbf}{\mathbf{Q}}
\nc{\Pbf}{\mathbf{P}}
\nc{\Stwo}{\mathbb{S}^2}
\nc{\tld}[1]{\tilde{#1}}
\nc{\wtld}[1]{\widetilde{#1}}
\nc{\hu}{\hat{u}}
\nc{\wh}[1]{\widehat{#1}}
\nc{\ph}{\varphi}
\nc{\sumeven}{\sum_{k=-N/2}^{N/2}{\hspace{-0.3cm}}'{\;\,}}
\nc{\sumevenk}{\sum_{k=-\frac{N_x}{2}}^{\frac{N_x}{2}}{\hspace{-0.3cm}}'{\;\,}}
\nc{\sumevenl}{\sum_{l=-\frac{N_y}{2}}^{\frac{N_y}{2}}{\hspace{-0.3cm}}'{\;\,}}
\nc{\sumevenm}{\sum_{m=-\frac{N_z}{2}}^{\frac{N_z}{2}}{\hspace{-0.3cm}}'{\;\,}}
\nc{\sumodd}{\sum_{k=-\frac{N-1}{2}}^{\frac{N-1}{2}}}
\nc{\sumoddl}{\sum_{l=-\frac{N-1}{2}}^{\frac{N-1}{2}}}
\nc{\cqfd}{~\hbox{\vrule width 2.5pt depth 2.5 pt height 3.5 pt}}
\nc{\rr}[1]{\textcolor{red}{#1}}
\nc{\ra}[1]{\renewcommand{\arraystretch}{#1}}
\nc{\nf}{\normalfont}
\title{Solving periodic semilinear stiff PDEs in 1D, 2D and 3D with exponential integrators}
\author{Hadrien Montanelli\thanks{Oxford University Mathematical Institute, Oxford OX2 6GG, UK.\ \ Supported by 
the European Research Council under the European Union's Seventh Framework Programme (FP7/2007--2013)/ERC grant agreement 
no.\ 291068.\ \ The views expressed in this article are not those of the ERC or the European Commission, and the European Union is not 
liable for any use that may be made of the information contained here.} 
\and Niall Bootland\thanks{Oxford University Mathematical Institute, Oxford OX2 6GG, UK.\ \ This publication was based on work supported in part by award no.\ KUK-C1-013-04,\ made by King Abdullah University of Science and Technology (KAUST).}}

\begin{document}

\maketitle

\begin{abstract}
Dozens of exponential integration formulas have been proposed for the high-accuracy solution of stiff PDEs such as 
the Allen--Cahn, Korteweg--de Vries and Ginzburg--Landau equations.
We report the results of extensive comparisons in MATLAB and Chebfun of such formulas in 1D, 2D and 3D, focusing on fourth and higher order methods,
and periodic semilinear stiff PDEs with constant coefficients.
Our conclusion is that it is hard to do much better than one of the simplest of these formulas, the ETDRK4 scheme of Cox and Matthews.
\end{abstract}



\pagestyle{myheadings}
\thispagestyle{plain}

\markboth{MONTANELLI AND BOOTLAND}{EXPONENTIAL INTEGRATORS FOR PERIODIC STIFF PDES}

\section{Introduction}

\begin{figure}
\includegraphics [scale=.5]{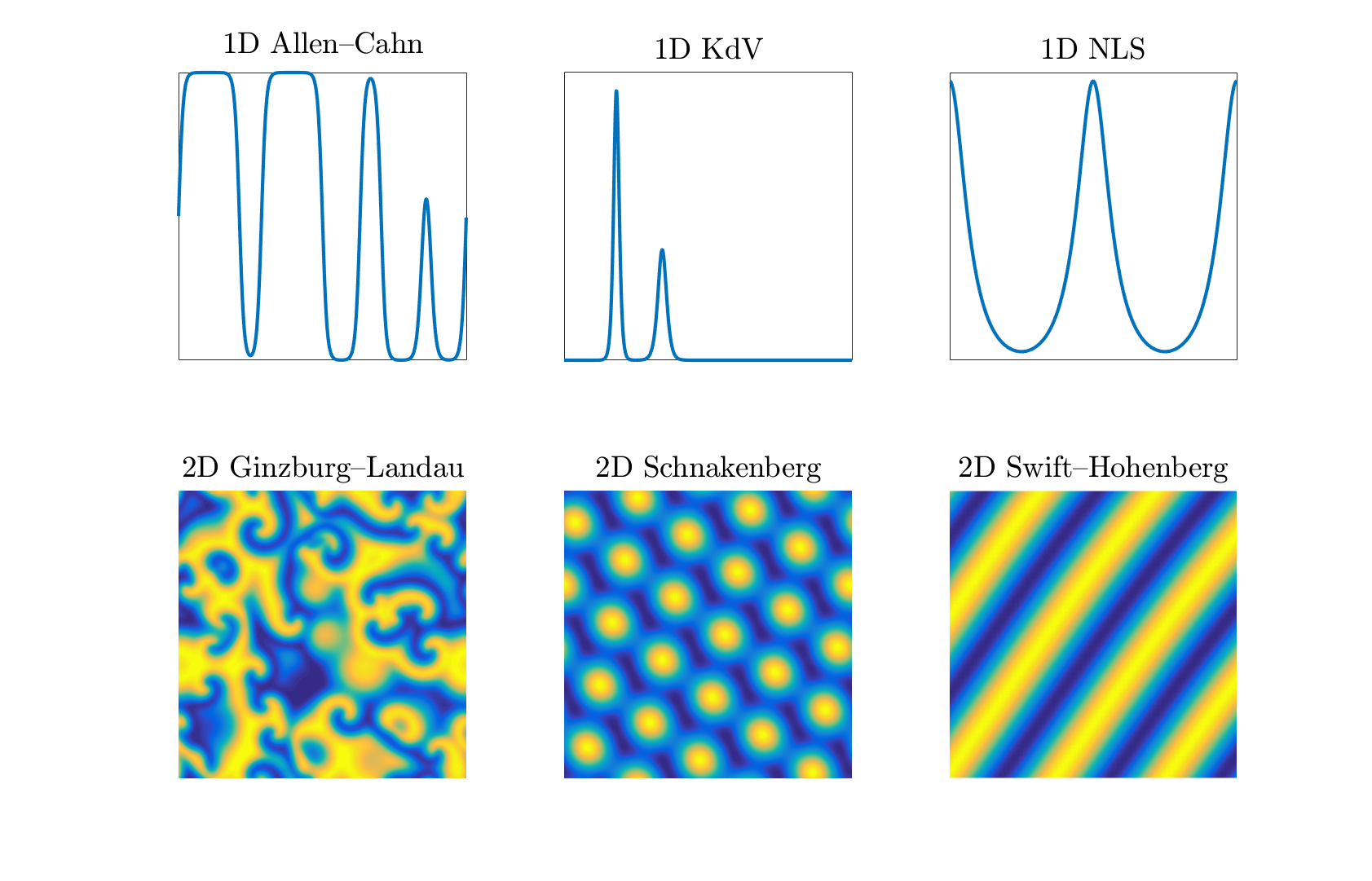}
\caption{\textit{First row (left to right): metastable solution of the Allen--Cahn equation~$\reff{AC}$, two-soliton solution of the KdV equation~$\reff{KdV}$ and breather solution of the NLS equation~$\reff{NLS}$.
Second row (left to right): frozen state solution of the Ginzburg--Landau equation~$\reff{GL}$, spot solution of the Schnakenberg equations~$\reff{Schnak}$ and convection rolls of the Swift--Hohenberg equations~$\reff{SH}$.}}
\label{iconic}
\end{figure}

We are interested in computing smooth solutions of stiff PDEs of the form
\begin{equation}
u_t = \mathcal{S}(u) = \mathcal{L}u + \mathcal{N}(u), \quad u(0,X)=u_0(X), \quad \text{periodic boundary conditions,}
\label{PDE}
\end{equation}

\noindent where $u(t,X)$ is a function of time $t$ and space $X$, $\mathcal{L}$ is a linear differential operator with constant coefficients on a domain in one, 
two or three space dimensions and $\mathcal{N}$ is a nonlinear differential (or non-differential) operator of lower order with constant coefficients and on the same domain.\footnote{$X$ 
denotes a space variable in 1D, 2D or 3D. Throughout this paper, we will use the variables $x$ in 1D, $(x,y)$ in 2D and $(x,y,z)$ in 3D. The domain will
be an interval in 1D, a square in 2D and a cube in 3D.}
In applications, PDEs of this kind typically arise when two or more different physical processes are combined, and many PDEs of interest in 
science and engineering take this form.  
For example, the Korteweg--de Vries equation $u_t = -u_{xxx} - uu_x$, the starting point of the study of nonlinear waves and solitons, 
couples third-order linear dispersion with first-order convection, and the Allen--Cahn equation $u_t = \epsilon u_{xx} + u - u^3$ couples second-order linear 
diffusion with a nondifferentiated cubic reaction term.  
Often a system of equations rather than a single scalar equation is involved, for example in the Gray--Scott and Schnakenberg equations, which 
involve two components coupled together.  
(The importance of coupling of nonequal diffusion constants in science was made famous by Alan Turing in the most highly-cited of all 
his papers~\cite{turing1952}.)
Fourth-order terms also arise, for example in the Cahn--Hilliard  equation, whose solutions describe structures of alloys, and in the 
Kuramoto--Sivashinsky equation, related to combustion problems among others, whose solutions are chaotic. 
Other examples of stiff PDEs include the Ginzburg--Landau, nonlinear Schr\"odinger (NLS) and Swift--Hohenberg equations.
Figure~\ref{iconic} shows six examples of solutions of such PDEs. 

Solving all these PDEs by generic numerical methods can be highly challenging. 
One of the main issues is \textit{stiffness}, characterised by the need for an explicit method to use small time-steps, 
much smaller than the condition required by accuracy. 
When too many steps are required, this can result in an infeasibly long computation.
There are other important issues such as the numerical conservation of various properties (e.g., for the KdV and NLS equations)---we shall not discuss this here; see, e.g.,~\cite{cano2015b}.

This paper describes and compares specialized methods that take advantage of two special features of~\reff{PDE}.
The first one is the periodic boundary conditions.
This allows us to discretize the spatial component of \reff{PDE} with a Fourier spectral method on $N$ points; equation~\reff{PDE} 
becomes a system of $N$ ODEs,
\begin{equation} 
\hat{u}' = \Sbf(\hat{u}) = \mathbf{L}\hat{u} + \mathbf{N}(\hat{u}), \quad \hat{u}(0)=\hat{u}_0,
\label{ODE}
\end{equation}

\noindent where $\hat{u}(t)$ is the vector of $N$ Fourier coefficients of the trigonometric interpolant of $u(t,X)$ at time $t$, 
and $\mathbf{L}$ (a $N\times N$ matrix) and $\mathbf{N}$ are the discretized versions of $\mathcal{L}$ and $\mathcal{N}$ in Fourier space. 
For example, in 1D on $[0, 2\pi]$ with $\mathcal{L}u=u_{xx}$ and an even number $N$ of equispaced grid points $\{x_j=2\pi j/N\}_{j=0}^{N-1}$,
we look for a solution $u(t,x)$ of the form\footnote{The prime on the summation sign in \reff{FourierSeries1D} signifies that the terms $k=\pm N/2$ are halved.}
\begin{equation}
u(t,x) \approx \sumeven \hat{u}_k(t) e^{ikx}
\label{FourierSeries1D}
\end{equation}

\noindent with Fourier coefficients 
\begin{equation} 
\hat{u}_k(t) = \frac{1}{N}\sum_{j=0}^{N-1}u(t,x_j)e^{-ikx_j}, \quad -\frac{N}{2}\leq k\leq \frac{N}{2}-1, \quad \hat{u}_{N/2}(t) = \hat{u}_{-N/2}(t).
\label{FourierCoeffs1D}
\end{equation}

\noindent Since FFT codes only store $N$ coefficients, the vector $\hat{u}(t)$ is defined as
\begin{equation}
\hat{u}(t) = \Big(\frac{\hat{u}_{-N/2}}{2}+\frac{\hat{u}_{N/2}}{2}, \hat{u}_{-N/2+1}(t),\ldots,\hat{u}_{N/2-1}(t)\Big)^T.
\label{vec2}
\end{equation}

\noindent For this PDE, $\Lbf=\mathbf{D}_N^{(2)}$ is the (diagonal) \textit{second-order Fourier differentiation matrix} with entries $-k^2$, $-N/2\leq k\leq N/2-1$.
In Section~3, we will also consider the first-, third-, and fourth-order Fourier differentiation matrices $\mathbf{D}_N$, $\mathbf{D}_N^{(3)}$, and $\mathbf{D}_N^{(4)}$;
see \cite{trefethen2000} for more details about Fourier spectral methods and \cite{montanelli2015b} for a review of trigonometric interpolation techniques.
(Note that stiffness is related to $\mathbf{L}$ having large eigenvalues since stability of spectral methods for time-dependent
PDEs requires that the eigenvalues of $\mathbf{L}$, scaled by the time-step, lie in the stability region of the time-stepping 
formula~\cite[Chapter 10]{trefethen2000}.) 

The second special feature of~\reff{PDE} is that it is semilinear, i.e., the higher-order terms of the equation are linear. 
Exponential integrators are a class of numerical methods for systems of ODEs that are aimed at taking advantage of this.
The linear part $\Lbf$, responsible for the stiffness, is integrated exactly using the matrix exponential while a numerical scheme is applied 
to $\mathbf{N}$. 

According to the 2005 review of Minchev and Wright~\cite{minchev2005}, the first exponential integrators were constructed by 
Certaine in 1960~\cite{certaine1960} and Pope in 1963~\cite{pope1963}.
Subsequently, however, Hochbruck and Ostermann~\cite{hochbruck2010} noted, in a comprehensive theoretical review of these schemes, that 
Hersch~\cite{hersch1958} had previously considered exponential integrators in 1958 in an effort to find schemes that are exact for linear problems with 
constant coefficients. 
The first use of the term \textit{exponential integrator} was by Hochbruck, Lubich and Selhofer~\cite{hochbruck1998} in a seminal paper of 1998.
The extensive use of these formulas for solving stiff PDEs seems to have been initiated by the papers by Cox and Matthews~\cite{cox2002} and 
Kassam and Trefethen~\cite{kassam2005}. 
A striking unpublished paper by Kassam~\cite{kassam2003} shows how effective such methods can be also for PDEs in 2D and 3D. 
A software package for such computations called EXPINT was produced by Berland, Skaflestad and Wright~\cite{berland2007}.

One of the simplest exponential integrators, commonly known as the Exponential Time Differencing (ETD)
Euler method, is given by\footnote{Throughout this paper, when introducing an exponential integrator such as \reff{ETDEuler}, $\hat{u}^n$ will mean $\hat{u}(t)$ at $t=t_n$.}
\begin{equation} 
\hat{u}^{n+1} = e^{h\Lbf}\hat{u}^n + h\ph_1(h\Lbf)\Nbf(\hat{u}^n),
\label{ETDEuler}
\end{equation}

\noindent where $h=t_{n+1}-t_n$ is the time-step and
\begin{equation} 
\ph_1(z) = \frac{e^z-1}{z}.
\label{phi1}
\end{equation}

\noindent As Minchev and Wright \cite{minchev2005} point out, this method has been rediscovered from many different viewpoints and has been known
by several other names. 
It can be derived by considering the linearized version of \reff{ODE} on $[t_n,t_{n+1}]$,
\begin{equation} 
\hat{u}' = \Sbf(\hat{u}^n) + \Sbf_{\hat{u}}(\hat{u}^n)(\hat{u} - \hat{u}^n), \quad \hat{u}(t_n)=\hat{u}^n,
\label{linODE}
\end{equation}

\noindent with exact solution at $t_{n+1} = t_n + h$,
\begin{equation} 
\hat{u}^{n+1} = \hat{u}^n + h\ph_1(h\Sbf_{\hat{u}}(\hat{u}^n))\Sbf(\hat{u}^n).
\label{ExpEuler}
\end{equation}

\noindent Approximating $\Sbf_{\hat{u}}(\hat{u}^n)$ by $\Lbf$ in \reff{ExpEuler} leads to \reff{ETDEuler}.
Note that \reff{ExpEuler} defines a time-stepping scheme too, known as the exponential Euler method.
The problem with \reff{ExpEuler} is that the exact Jacobian $\Sbf_{\hat{u}}$, and its value under the exponential-like function $\ph_1$, need to be 
computed at each time-step.
This would involve a high computational cost, so typically one either does not compute $\ph_1$ but rather an approximation,
such as a Pad\'e approximation, or else one uses an approximation to the Jacobian as opposed to the exact Jacobian.
The former approach contains the Rosenbrock methods~\cite{hairer1982, hairer1991, hochbruck2009, luan2014b, rosenbrock1963, vanderhouwen1977} 
and the Exponential Propagation Iterative methods of Runge-Kutta type (EPIRK)~\cite{rainwater2014, rainwater2015, tokman2006, tokman2011, tokman2012a}.
The latter approach is what we shall consider in this paper.

As we just described, exponential integrators are characterised by the use of exponential and related functions of the matrix $\mathbf{L}$. 
Standard methods for computing the matrix exponential in the context of exponential integrators include the scaling and squaring method~\cite{almohy2011}, 
the Carath\'{e}odory--Fej\'{e}r method~\cite{schmelzer2007} and Krylov subspace methods~\cite{tokman2006}.
There is recent work that shows that exponential integrators together with Krylov methods are competitive, for instance see Tokman and Loffeld~\cite{loffeld2013, tokman2012b}.
In our case we consider periodic problems with constant coefficients so the matrices are diagonal and the matrix exponential is trivial.

We compare in this paper 30 exponential integrators of fourth and higher order on 11 model problems in 1D, 2D and 3D, 
using MATLAB R2015b and Chebfun v5.5~\cite{chebfun}.
Comparisons with other types of time-stepping schemes are out of the scope of the article; see, e.g.,~\cite{garcia2014, kassam2003, klein2008, klein2011, loffeld2013}.
Let us emphasize that we are interested in determining if one of the high order integrators outperforms the others on a large class of problems.
For a particular problem, it might be possible to design a very specific scheme, of possibly lower order than four, which performs extremely well. 
For example, Cano and Gonz\'{a}les--Pach\'{o}n have recently shown that the low-order Lawson methods, combined with orthogonal projections onto 
some invariants, can be very competitive for the nonlinear Schr\"{o}dinger equation~\cite{cano2015a, cano2015b}.
Let us also emphasize that since we only consider periodic problems, we do not expect to see any order reduction in the convergence of the exponential integrators, as already
observed in, e.g., \cite{kassam2005}. 
For different types of boundary conditions (e.g., homogeneous Dirichlet conditions), certain schemes (e.g., Lawson methods) do not satisfy 
the so-called \textit{stiff order conditions}~\cite{hochbruck2010}---which guarantee a certain order of convergence independently of the considered problem---and can therefore exhibit a strong order reduction in practice.

The paper is structured as follows.
We present the 30 exponential integrators in Section~2 and the 11 model problems in Section~3.
The numerical results are presented in Section~4 and show that it is hard to do much better than one of the simplest of these formulas, the ETDRK4 scheme of Cox and 
Matthews~\cite{cox2002}.

\section{Thirty exponential integrators}

\subsection{Exponential general linear methods}

We consider exponential integrators, based on the approximation of the Jacobian of \reff{ExpEuler}, that belong to the large class 
of exponential general linear methods, first introduced by Minchev and Wright in 2005~\cite{minchev2005}.
This class contains, in particular, the ETD Runge--Kutta (one-step), ETD Adams--Bashforth (multistep), Lawson and exponential predictor-corrector methods.
For given starting values $\hat{u}^0,\hat{u}^1,\ldots,\hat{u}^{q-1}$ at times $t=0,h,\ldots,(q-1)h$, the numerical approximation $\hat{u}^{n+1}$ at time 
$t_{n+1}=(n+1)h$, $n+1\geq q$, is given by the formula
\begin{equation}
\dsp \hat{u}^{n+1} = e^{h\Lbf}\hat{u}^n + h\sum_{i=1}^{s}B_i(h\Lbf) \Nbf(\hat{v}^i) +  h\sum_{i=1}^{q-1}V_i(h\Lbf) \Nbf(\hat{u}^{n-i}),
\label{step1}
\end{equation}

\noindent with $q$ \textit{steps} $\hat{u}^{n-i}$ and $s$ \textit{stages} $\hat{v}^i$, with $\hat{v}^1=\hat{u}^n$ and
\begin{equation}
\dsp \hat{v}^i = e^{C_ih\Lbf}\hat{u}^n + h\sum_{j=1}^{i-1}A_{i,j}(h\Lbf) \Nbf(\hat{v}^j) + h\sum_{j=1}^{q-1}U_{i,j}(h\Lbf) \Nbf(\hat{u}^{n-j}), \quad 2\leq i\leq s.
\label{step2}
\end{equation}

\begin{table}
\caption{\textit{Butcher tableau of an exponential integrator with $q$ steps and $s$ stages.}}
\vspace{.5cm}
\centering
\ra{1.3}
\begin{tabular}{c|cccc|ccc}
\toprule
$C_2$     & $A_{2,1}$  &                   &                   &                   & $U_{2,1}$ &  $\ldots$ & $U_{2,q-1}$ \\
$\vdots$ & $\vdots$    & $\ddots$    &                   &                   & $\vdots$   &               &  $\vdots$      \\
$C_s$     & $A_{s,1}$   & $\ldots$     & $A_{s,s-1}$ &                  & $U_{s,1}$  & $\dots$  & $U_{s,q-1}$   \\            
\midrule
              & $B_1$         & $\ldots$     & $B_{s-1}$ & $B_s$          & $V_1$       & $\dots$  & $V_{q-1}$ \\
\bottomrule
\end{tabular}
\label{Butcher}
\end{table}

\noindent Each scheme is characterised by its coefficients $A$, $B$, $C$, $U$, and $V$, which can be conveniently listed in a Butcher tableau, 
as in Table~\ref{Butcher}. Note that these coefficients (except $C$) depend on $\Lbf$---for instance, \reff{ETDEuler} uses one stage and one step,
and its only non-zero coefficient is $B_1=\ph_1(h\Lbf)$.
Note that, in practice, the nonlinear evaluations $\Nbf(\hat{v}^i)$ and $\Nbf(\hat{u}^{n-i})$ are carried out in value space, e.g.,
$\Nbf(\hat{v}^i)$ means $\mathbf{F}\big(\Nbf\big(\mathbf{F}^{-1}\hat{v}^i\big)\big)$, with discrete Fourier transform $\mathbf{F}$.
Methods of the form \reff{step1}--\reff{step2} not only include purely one-step methods 
($q=1$, $s\geq1$) and purely multistep methods ($q\geq1$, $s=1$), but also combinations of both.

The coefficients satisfy the following summation properties,\footnote{There are two exceptions: the coefficients of 
the Lawson4 and ABLawson4 schemes do not satisfy the summation properties~\reff{Summation}.}
\begin{equation}
\begin{array}{l}
\dsp B_1 = \ph_1(h\Lbf) - \sum_{i=2}^{s}B_i(h\Lbf) - \sum_{i=1}^{q-1}V_i(h\Lbf), \\\\
\dsp A_{i,1} = \psi_{1,i}(h\Lbf) - \sum_{j=2}^{i-1}A_{i,j}(h\Lbf) - \sum_{j=1}^{q-1}U_{i,j}(h\Lbf), \quad 2\leq i\leq s,
\end{array}
\label{Summation}
\end{equation}

\noindent where the $\ph$- and $\psi$-functions are exponential and related functions that we shall define in the next subsection.
As a consequence, it is notationally convenient to incorporate this condition by filling the corresponding entries of the Butcher tableau with a dot on the understanding that these method coefficients are given by \reff{Summation}.
Note that exponential integrators of the form \reff{step1}--\reff{step2} do not include the EMAM4 scheme of Calvo and Palencia~\cite{calvo2006}. 
It has been shown in~\cite{bootland2014} that it often suffers from stability problems.

Let us finish this section with a few words about the computational cost per time-step.
Since the matrices in \reff{step1}--\reff{step2} are diagonal, the matrix-vector products cost only $\mathcal{O}(N)$ operations.
The dominant cost per time-step is then the cost of an FFT, i.e., $\mathcal{O}(N\log N)$ operations.
For exponential integrators of the form~\reff{step1}--\reff{step2}, the total cost to compute $\hat{u}^{n+1}$ is therefore $\mathcal{O}(2sN\log N)$.
As a consequence, purely multistep methods have a low computational cost per time-step.

\subsection{Evaluating the $\ph$-functions}

The coefficients $A$, $B$, $C$, $U$, and $V$ involve the $\ph$ and $\psi$-functions applied to $\Lbf$. 
Because $\Lbf$ is diagonal, $\ph(\Lbf)$ and $\psi(\Lbf)$ reduce to $\ph$ and $\psi$ applied to the diagonal elements $\lambda$ of $\Lbf$, 
so all we have to be able to is to compute $\ph(\lambda)$ and $\psi(\lambda)$ for $\lambda\in\C$.
The $\ph$-functions are defined by the recurrence relation,
\begin{equation}
\ph_{l+1}(z) = \frac{\ph_l(z) - 1/l!}{z}, \quad l\geq 1,
\label{phi}
\end{equation}

\noindent with $\ph_0(z)=e^z$. After $\ph_0$, the first few $\ph$-functions are \reff{phi1} and
\begin{equation}
\ph_2(z) = \frac{e^z - z - 1}{z^2}, \quad \ph_3(z) = \frac{e^z - \frac{z^2}{2} - z - 1}{z^3},
\end{equation}

\noindent while the $\psi$-functions are defined via the $\ph$-functions and the coefficients $C$,
\begin{equation}
\psi_{l,m}(z) = C_m^l\ph_l(C_m z), \quad l\geq 0, \quad 1\leq m\leq s.
\label{psi}
\end{equation}

\noindent Equations \reff{phi} and \reff{psi} can be implemented recursively, but the accurate evaluation of $\ph$ and $\psi$
is not straightforward because it can suffer from cancellation error.
Following the idea of Kassam and Trefethen in~\cite{kassam2005}, to compute the functions at some $\lambda\in\C$, we use 
Cauchy's integral formula
\begin{equation}
\ph(\lambda) = \frac{1}{2\pi i}\oint_\Gamma \frac{\ph(z)}{z - \lambda} dz,
\end{equation}

\noindent which can be approximated with exponential accuracy by the trapezoidal rule \cite{trefethen2014},
\begin{equation}
\ph(\lambda) \approx \frac{1}{M} \sum_{k=1}^M \ph\big(\lambda + e^{2\pi i (k-0.5)/M}\big),
\label{evalphi1}
\end{equation}

\noindent taking $\Gamma$ to be the circle of radius $1$ centred at $\lambda$, oriented counter-clockwise, discretized with $M$ equally spaced points. 
Note that the $\ph$-functions satisfy $\ph(\bar{z})=\bar{\ph}(z)$ for all $z\in\C$. 
As a consequence, when $\lambda$ is on the real axis, we can take $\Gamma$ to be the upper half of the circle of radius~$1$ centred at $\lambda$
and take the real part of the result, i.e., 
\begin{equation}
\dsp \ph(\lambda) = \frac{1}{\pi}\mathcal{R}\Bigg(\int_0^{\pi} \ph\big(\lambda + e^{i\theta}\big) d\theta\Bigg),
\end{equation}

\noindent which can be discretized by
\begin{equation}
\ph(\lambda) \approx \frac{1}{M} \mathcal{R}\Bigg(\sum_{k=1}^M \ph\big(\lambda + e^{\pi i (k-0.5)/M}\big)\Bigg).
\label{evalphi2}
\end{equation}

\noindent If this symmetry is not explicitly used in the computation of the $\ph$-functions when $\lambda$ is real,
rounding errors appear that lead to numerical instability. 
Note that the evaluation of the $\ph$-functions using \reff{evalphi1} or \reff{evalphi2} requires $\mathcal{O}(M)$ operations per $\lambda$.

Let us emphasize that we can use circles of radius $1$ around each eigenvalue $\lambda$ of $\Lbf$ because $\Lbf$ is diagonal. 
When the matrix is not diagonal, one has to use a single contour that encloses all the eigenvalues, and the best possible contour depends on the problem. 
For example, for parabolic problems (also called diffusive problems), all the eigenvalues are on the real negative axis and the best contour 
is a Hankel contour~\cite{schmelzer2007, trefethen2006}.

Using contour integrals is not the only possible remedy for cancellation error. When Pope introduced~\reff{ETDEuler}, he suggested the
use of Taylor series for small~$\lambda$ and the direct formula for large~$\lambda$. 
The problem with this approach is that, for some intermediate values, neither method gives full precision, 
as noted by Cox and Matthews~\cite{cox2002} and Kassam and Trefethen~\cite{kassam2005}. 
Another approach is to use Pad\'e approximations, combined with a scaling and squaring technique~\cite{beyklin1998}. 
This method is also effective, but the contour integral method is particularly appealing because of its greater generality for dealing with arbitrary functions.

\subsection{Introducing the 30 exponential integrators}

\begin{table}
\caption{\textit{A reference table for the exponential integrators considered in this paper.
Since there is no order reduction for periodic (diagonal) problems, the stiff convergence order is the same at the non-stiff convergence order.
Note that some methods do not appear explicitly in the references listed but can be derived using order conditions or recurrence formulas found there; 
these cases are marked by asterisks.
The Butcher tableaux can be found in the Ph.D. thesis of the first author \cite{montanelli2017}.}}
\vspace{.5cm}
\centering
\ra{1.2}
\small
\begin{tabular}{cccccc}
\toprule
\textbf{Method} & \textbf{Type} & \textbf{Order} & \textbf{Stages} $s$ & \textbf{Steps} $q$ & \textbf{Ref.} \\
\midrule
ABN{\o}rsett4 & ETD Adams--Bashforth & 4 & 1 & 4 & \cite{norsett1969}$^{\ast}$ \\
ABN{\o}rsett5 & ETD Adams--Bashforth & 5 & 1 & 5 &  \cite{norsett1969}$^{\ast}$ \\
ABN{\o}rsett6 & ETD Adams--Bashforth & 6 & 1 & 6 &  \cite{norsett1969}$^{\ast}$ \\
\midrule
ETDRK4 & ETD Runge--Kutta & 4 & 4 & 1 & \cite{cox2002}\hphantom{$^{\ast}$} \\
Friedli & ETD Runge--Kutta & 4 & 4 & 1 & \cite{friedli1978}\hphantom{$^{\ast}$} \\
Krogstad & ETD Runge--Kutta & 4 & 4 & 1 & \cite{krogstad2005}\hphantom{$^{\ast}$} \\
Minchev & ETD Runge--Kutta & 4 & 4 & 1 & \cite{minchev2004}\hphantom{$^{\ast}$} \\
Strehmel--Weiner & ETD Runge--Kutta & 4 & 4 & 1 & \cite{strehmel1982}\hphantom{$^{\ast}$} \\
Hochbruck--Ostermann & ETD Runge--Kutta & 4 & 5 & 1 & \cite{hochbruck2005}\hphantom{$^{\ast}$} \\
EXPRK5S8 & ETD Runge--Kutta & 5 & 8 & 1 & \cite{luan2014a}\hphantom{$^{\ast}$} \\
\midrule
ABLawson4 & Lawson & 4 & 1 & 4 & \cite{lawson1967}\hphantom{$^{\ast}$} \\
Lawson4 & Lawson & 4 & 4 & 1 &  \cite{lawson1967}\hphantom{$^{\ast}$} \\
\midrule
GenLawson41 & Gen. Lawson & 4 & 4 & 1 & \cite{krogstad2005}\hphantom{$^{\ast}$} \\
GenLawson42 & Gen. Lawson & 4 & 4 & 2 & \cite{krogstad2005}\hphantom{$^{\ast}$} \\
GenLawson43 & Gen. Lawson & 4 & 4 & 3 & \cite{krogstad2005}\hphantom{$^{\ast}$} \\
GenLawson44 & Gen. Lawson & 5 & 4 & 4 & \cite{krogstad2005}\hphantom{$^{\ast}$} \\
GenLawson45 & Gen. Lawson & 6 & 4 & 5 & \cite{krogstad2005}\hphantom{$^{\ast}$} \\
\midrule
ModGenLawson41 & Mod. Gen. Lawson & 4 & 4 & 1 & \cite{ostermann2006}$^{\ast}$ \\
ModGenLawson42 & Mod. Gen. Lawson & 4 & 4 & 2 & \cite{ostermann2006}$^{\ast}$ \\
ModGenLawson43 & Mod. Gen. Lawson & 4 & 4 & 3 & \cite{ostermann2006}$^{\ast}$ \\
ModGenLawson44 & Mod. Gen. Lawson & 5 & 4 & 4 & \cite{ostermann2006}$^{\ast}$ \\
ModGenLawson45 & Mod. Gen. Lawson & 6 & 4 & 5 & \cite{ostermann2006}$^{\ast}$ \\
\midrule
PEC423      & Exp. Predictor-Corrector & 4 & 2 & 3 & \cite{ostermann2006}$^{\ast}$ \\
PECEC433 & Exp. Predictor-Corrector & 4 & 3 & 3 & \cite{ostermann2006}$^{\ast}$ \\
PEC524     & Exp. Predictor-Corrector & 5 & 2 & 4 & \cite{ostermann2006}$^{\ast}$ \\
PECEC534  & Exp. Predictor-Corrector & 5 & 3 & 4 & \cite{ostermann2006}$^{\ast}$ \\
PEC625      & Exp. Predictor-Corrector & 6 & 2 & 5 & \cite{ostermann2006}$^{\ast}$ \\
PECEC635  & Exp. Predictor-Corrector & 6 & 3 & 5 & \cite{ostermann2006}$^{\ast}$ \\
PEC726      & Exp. Predictor-Corrector & 7 & 2 & 6 & \cite{ostermann2006}$^{\ast}$ \\
PECEC736  & Exp. Predictor-Corrector & 7 & 3 & 6 & \cite{ostermann2006}$^{\ast}$ \\
\bottomrule
\end{tabular}
\label{ReferenceTable}
\end{table}

Table~\ref{ReferenceTable} lists the exponential integrators considered in this paper. 
Their Butcher tableaux can be found in the Ph.D. thesis of the first author \cite{montanelli2017}.

\paragraph{ETD Adams--Bashforth} The first category of exponential integrators is the \textit{ETD Adams--Bashforth} schemes of order four to six.
These are ETD (purely) multistep methods, which reduce to Adams--Bashforth schemes when $\Lbf=0$, 
and go back to N{\o}rsett in 1969~\cite{norsett1969}. 
Since $s=1$, \reff{step1}--\reff{step2} takes the simpler form
\begin{equation}
\dsp \hat{u}^{n+1} = e^{h\Lbf}\hat{u}^n + hB_1(h\Lbf) \Nbf(\hat{u}^n) +  h\sum_{i=1}^{q-1}V_i(h\Lbf) \Nbf(\hat{u}^{n-i}),
\label{ETDAB}
\end{equation}

\noindent i.e., the only non-zero coefficients are $B_1$ and those in $V$.
Note that since these schemes are purely multistep, \reff{ETDAB} only requires two FFTs per time-step.
We label these methods as ABN{\o}rsett$q$, where $4\leq q\leq 6$ is the order and also the number of steps.
For more details on the derivation of ABN{\o}rsett methods see Minchev and Wright~\cite{minchev2005}, 
who also show a connection between these schemes and the IMEX schemes of Ascher, Ruuth and Wetton~\cite{ascher1995}. 
One can also derive methods based on Adams--Moulton methods, known as AMN{\o}rsett$q$ methods.
These are implicit but can be used within predictor-corrector pairs, as we will see when introducing exponential predictor-corrector schemes (the last category in the table). 
A comprehensive look at both the Adams--Bashforth and Adams--Moulton exponential integrators can be found in the paper by Hochbruck and Ostermann~\cite{hochbruck2011}.

\begin{table}
\caption{\textit{Butcher tableau for ETDRK$\,4$. Note the dots in the first column, which indicate that these coefficients are computed using $\reff{Summation}$.
These coefficients are (from top to bottom): $A_{2,1}=\psi_{1,2}$, $A_{3,1}=0$, $A_{4,1}=\ph_1 - 2\psi_{1,2}$ and $B_1=\ph_1 - 3\ph_2 + 4\ph_3$.}}
\vspace{.5cm}
\centering
\ra{1.5}
\begin{tabular}{c|ccccc|c}    
\toprule
$\frac{1}{2}$ & $A_{2,1}=\cdot$ & & & & & \\
$\frac{1}{2}$ & $A_{3,1}=\cdot$ & $A_{3,2}=\psi_{1,2}$ & & & & \\
$1$ & $A_{4,1}=\cdot$ & $A_{4,2}=0$ & $A_{4,3}=2\psi_{1,2}$ & & & \\
\midrule
& $B_1=\cdot$ & $B_2=2\ph_2 - 4\ph_3$ & $B_3=2\ph_2 - 4\ph_3$ & $B_4=-\ph_2 + 4\ph_3$ & & \\
\bottomrule
\end{tabular}
\label{ETDRK4tableau}
\end{table}

\paragraph{ETD Runge--Kutta} The second category is the \textit{ETD Runge--Kutta} schemes of order four to five.
These are (purely) one-step methods and go back to Friedli in 1978~\cite{friedli1978} and Strehmel--Weiner in 1982~\cite{strehmel1982}.
More recently, inspired by Cox and Matthews' ETDRK4 scheme~\cite{cox2002}, Minchev~\cite{minchev2004} in 2004 and
Krogstad~\cite{krogstad2005} and Hochbruck and Ostermann~\cite{hochbruck2005} in 2005 derived ETD Runge--Kutta schemes of order four.
Luan and Ostermann proposed a scheme of order five (EXPRK5S8) in 2014~\cite{luan2014a}.
Overviews of ETD Runge--Kutta methods and some of their history can be found in the reviews of Hochbruck and Ostermann~\cite{hochbruck2010} and Minchev and Wright~\cite{minchev2005}, where connections are described between ETD Runge--Kutta, generalised Runge--Kutta and semi-implicit methods.
Since $q=1$, \reff{step1}--\reff{step2} reduces to
\begin{equation}
\begin{array}{l}
\dsp \hat{u}^{n+1} = e^{h\Lbf}\hat{u}^n + h\sum_{i=1}^{s}B_i(h\Lbf) \Nbf(\hat{v}^i), \\
\dsp \hat{v}^1=\hat{u}^n, \quad \hat{v}^i = e^{C_ih\Lbf}\hat{u}^n + h\sum_{j=1}^{i-1}A_{i,j}(h\Lbf) \Nbf(\hat{v}^j), \quad 2\leq i\leq s.
\end{array}
\label{ETDRK}
\end{equation}

\noindent The only non-zero coefficients are those in $A$, $B$ and $C$. 
The coefficients for the ETDRK4 scheme can be found in Table~\ref{ETDRK4tableau} and correspond to the following formula:
\begin{equation}
\begin{array}{l}
\hat{v}^1 = \hat{u}^n, \\\\
\hat{v}^2 = e^{\Lbf h/2}\hat{u}^n + (h/2)\ph_1(\Lbf h/2)\Nbf(\hat{v}^1), \\\\
\hat{v}^3 = e^{\Lbf h/2}\hat{u}^n + (h/2)\ph_1(\Lbf h/2)\Nbf(\hat{v}^2), \\\\
\hat{v}^4 = e^{\Lbf h/2}\hat{v}^2 + (h/2)\ph_1(\Lbf h/2)[2\Nbf(\hat{v}^3) - \Nbf(\hat{v}^1)], \\\\
\hat{u}^{n+1} = e^{\Lbf h}\hat{u}^n + hB_1\Nbf(\hat{v}^1) + hB_2[\Nbf(\hat{v}^2) + \Nbf(\hat{v}^3)] + hB_4\Nbf(\hat{v}^4),
\end{array}
\label{ETDRK4}
\end{equation}

\noindent where $B_1=B_1(h\Lbf),\,\ldots,B_4=B_4(h\Lbf)$.

\paragraph{Lawson} The third category is the \textit{Lawson} methods. 
First developed by Lawson in 1967~\cite{lawson1967}, and often known as integrating factor (IF) methods,
the motivation behind Lawson methods is to use a change of variable in~\reff{ODE} to get rid of the stiff linear part, 
and then apply a numerical solver to the transformed equation. 
The Lawson transformation consists of the change of variables $\hat{v}(t)=e^{-\Lbf t}\hat{u}(t)$.
If we differentiate this and substitute into~\reff{ODE}, the transformed equation is
\begin{equation}
\hat{v}' = e^{-\Lbf t}\Nbf(e^{\Lbf t}\hat{v}), \quad \hat{v}(0)=\hat{u}_0.
\label{ODELawson}
\end{equation}

\noindent The linear term is gone, and the transformed equation \reff{ODELawson}, while no longer stiff, now has rapidly varying coefficients.
Once we have decided on a scheme to solve \reff{ODELawson}, we can transform back to $\hat{u}$. Lawson, in his 1967 paper, used the classical fourth order Runge--Kutta scheme
on the transformed equation \reff{ODELawson}; we call this method Lawson4. 
Using the classical fourth order Adams--Bashforth scheme gives the ABLawson4 method.
Ehle and Lawson observed in \cite{ehle1975} that Runge--Kutta based Lawson methods only work well when the problem is moderately stiff.
Another problem with Lawson methods, as indicated by Krogstad~\cite{krogstad2005}, is that they do not preserve fixed points of the differential equation. 

\paragraph{Generalised Lawson} Krogstad worked around these problems to derive \textit{generalised Lawson methods,} also called generalised integrating factor (GIF) methods, the fourth category in the table.
These are based on the transform
\begin{equation}
\hat{v}(t) = e^{-\Lbf t}\hat{u}(t) - e^{-\Lbf t}\sum_{l=1}^{q}t^{l}\ph_{l}(t\Lbf)p_{l-1}, 
\end{equation}

\noindent where the $p_{l}$ are the coefficients, in a (scaled) monomial basis, of the polynomial $P(t)$ of degree $q-1$ that interpolates 
the values $\{N(\hat{u}^{n-l})\}_{l=1}^q$ at the points $\{t_{n-l}\}_{l=1}^q$; see \cite{hochbruck2010, minchev2005} for details. 
Differentiating this and substituting into~\reff{ODE} leads to the transformed equation
\begin{equation}
\hat{v}' = e^{-\Lbf t}\Big(\Nbf\Big(e^{\Lbf t}\hat{v} + \sum_{l=1}^{q}t^{l}\ph_{l}(t\Lbf)p_{l-1}\Big)-P(t)\Big), \quad \hat{v}(0)=\hat{u}_0.
\label{ODEGenLawson}
\end{equation}

\noindent Note that \reff{ODELawson} is the special case of \reff{ODEGenLawson} with $P(t)=0$.
The idea of Krogstad is to apply, for various values of $q$, the classical fourth order Runge--Kutta scheme on \reff{ODEGenLawson}, and  then transform
back to $\hat{u}$.
It leads to methods with four stages and $q$ steps, called the GenLawson$4q$ methods.

\paragraph{Modified generalised Lawson} As we increase $q$ in the generalised Lawson methods we incorporate more of the nonlinearity and the methods have improved accuracy. 
However, this in part comes at the cost of stability, especially for dispersive problems, as was demonstrated by Krogstad~\cite{krogstad2005}. 
A modification, based on satisfying order conditions, given by Ostermann, Thalhammer and Wright~\cite{ostermann2006}, significantly improves stability. 
The modification is given by the requirement that
\begin{equation}
\sum_{i=1}^{4} B_{i}(h\Lbf)\frac{c_{i}^{j}}{j!} + \sum_{i=1}^{q-1} V_{i}(h\Lbf)\frac{(-1)^{j}}{j!} = \ph_{j+1}(h\Lbf), \quad 0\leq j\leq q-1,
\end{equation}
where, as before, $q-1$ is the degree of the polynomial approximation. These are the \textit{modified generalised Lawson} methods, labelled as ModGenLawson$4q$.

\paragraph{Exponential predictor-corrector} Just as with the standard Adams--Bashforth and Adams--Moulton multistep methods, the exponential versions can be used in predictor-corrector pairs. 
These are the \textit{exponential predictor-corrector} methods, the last category in the table.
For instance, using ABN{\o}rsett3 for a predictor step and AMN{\o}rsett4 for the corrector step yields the fourth order method called PEC423 in the MATLAB package EXPINT~\cite{berland2007}.
(PEC stands for predict-evaluate-correct, four is the order, two is the number of stages and three is the number of steps.)
One can evaluate and correct again, that is, use the corrector twice. 
The name PECEC433 is given in EXPINT for the fourth order method that uses ABN{\o}rsett3 for a predictor step and AMN{\o}rsett4 for two corrector steps. 

\section{Eleven model problems}

In this section we describe the PDEs used in the comparisons of Section 4, including the initial conditions, the domains and the time intervals.
There are five PDEs in 1D and three PDEs considered in both 2D and 3D; see Table~\ref{PDEs}.

\begin{table}
\caption{\textit{The model problems we consider in this paper. The linear operator of a diffusive PDE has real eigenvalues while
it has purely imaginary eigenvalues for dispersive PDEs. Note that we take $A=0$ for the Ginzburg--Landau equation~$\reff{GL}$. 
For $A\neq0,$ the linear part would have complex eigenvalues.}}
\vspace{.5cm}
\centering
\ra{1.3}
\begin{tabular}{ccc}
\toprule
PDE & Dimension & Stiff linear Part\\
\midrule
Allen--Cahn & 1D & second-order diffusive \\
Cahn--Hilliard & 1D & fourth-order diffusive \\
Korteweg--de Vries & 1D & third-order dispersive \\
Kuramoto--Sivashinsky & 1D & fourth-order diffusive \\
nonlinear Schr\"odinger & 1D  & second-order dispersive \\
Ginzburg--Landau & 2D \& 3D & second-order diffusive \\
Schnakenberg & 2D \& 3D & second-order diffusive \\
Swift--Hohenberg & 2D \& 3D & fourth-order diffusive \\
\bottomrule
\end{tabular}
\label{PDEs}
\end{table}

\subsection{Model problems in 1D}

\paragraph{Allen--Cahn} The \textit{Allen--Cahn equation}, derived by Allen and Cahn in the 1970s, is a reaction-diffusion equation which
describes the process of phase separation in iron alloys (see, e.g., \cite{allen1979}). It is given in one dimension as
\begin{equation}
u_t = \epsilon u_{xx} + u - u^3,
\label{AC}
\end{equation}

\noindent with linear diffusion $\epsilon u_{xx}$ and a cubic reaction term $u-u^3$.
The function $u$ is the order parameter, a correlation function related to the positions of the different components of the alloy.
The Allen--Cahn equation exhibits stable equilibria at $u=\pm 1$ while $u=0$ is an unstable equilibrium. 
Solutions often display metastability where wells $u\approx -1$ compete with peaks $u\approx 1$,
and structures remain almost unchanged for long periods of time before changing suddenly.
This can be quantified: features with width $L$ persist for time scales on the order of $e^{L/\epsilon}$.
In Fourier space with a grid of size $N$, \reff{AC} becomes 
\begin{equation}
\hat{u}' = \epsilon \mathbf{D}^{(2)}_N\hat{u} + \hat{u} - \mathbf{F}\big(\big(\mathbf{F}^{-1}\hat{u}\big)^3\big).
\end{equation}

\noindent We take $\epsilon=5\times10^{-2}$, 
\begin{equation}
u(0,x) = \frac{1}{3}\tanh(2\sin(x)) - e^{-23.5.(x-\pi/2)^2} + e^{-27(x-4.2)^2} + e^{-38(x-5.4)^2}, 
\label{ACIC}
\end{equation}

\noindent with $x\in[0,2\pi]$ and solve up to $t=60$.
This initial condition quickly converges to a set of wells $u\approx -1$ and peaks $u\approx 1$ (at around $t=4$) and
eventually to a two-plateau solution (at around $t=500$).
Figure~\ref{iconic} shows the solution at time $t=113$, when the peak on the far right is switching to $u\approx-1$.

\paragraph{Cahn--Hilliard} The \textit{Cahn--Hilliard equation},
\begin{equation}
u_t = \alpha(-u_{xx} - \gamma u_{xxxx} + (u^3)_{xx}),
\label{CH}
\end{equation}

\noindent is a fourth order reaction-diffusion problem which Cahn and Hilliard proposed in 1958 as a model for the process of phase separation in binary alloys~\cite{cahn1958}.
It couples second-order destabilizing diffusion $-u_{xx}$ with fourth-order stabilizing diffusion $-u_{xxxx}$ and a differentiated cubic reaction term $(u^3)_{xx}$.
The function $u$ is defined as $u=1-2c_A$ where $0\leq c_A\leq1$ and $c_B=1-c_A$ denote the concentrations of the two components $A$ and $B$ of the alloy, that is,
$u=-1$ means pure $A$ while $u=1$ means pure $B$.
The Cahn--Hilliard equation also exhibits metastable solutions.
When quenched below a critical temperature, alloys described by \reff{CH} become unstable in the sense that small metastable pockets of relatively pure $A$ and $B$ may soon appear, corresponding to 
wells $u=-1$ and peaks $u=1$. These pockets may coarsen into larger pockets at progressively larger times. 
In Fourier space, \reff{CH} becomes
\begin{equation}
\hat{u}' = \alpha(-\Dbf_N^{(2)} - \gamma \Dbf_N^{(4)})\hat{u} + \alpha\Dbf_N^{(2)}\mathbf{F}\big(\big(\mathbf{F}^{-1}\hat{u}\big)^3\big).
\end{equation}

\noindent We take $\alpha=10^{-2}$, $\gamma=10^{-3}$,
\begin{equation}
u(0,x) = \frac{1}{5}\sin(4\pi x)^5 - \frac{4}{5}\sin(\pi x), \quad x\in[-1,1],
\end{equation}

\noindent and solve up to $t=12$. This initial condition evolves to a four-plateau solution (two wells $u\approx -1$, two peaks $u\approx 1$) at around $t=12$ before switching to a two-plateau solution (one well, one peak)
at around $t=70$.

\paragraph{Korteweg--de Vries} The \textit{KdV equation},
\begin{equation}
u_t = - u_{xxx} - uu_x,
\label{KdV}
\end{equation}

\noindent was derived by Korteweg and de Vries in 1895 to model the propagation of waves in shallow water \cite{korteweg1895}.
It couples dispersion $-u_{xxx}$ with nonlinear convection $-uu_x$. Among the solutions of \reff{KdV} are solitary waves or \textit{solitons}.
These are waves that maintain their shapes as they travel and are given by 
\begin{equation}
u(t,x) = 3c\,\mrm{sech}^2\Big(\frac{\sqrt{c}}{2}(x-x_0-ct)\Big), \quad c>0.
\label{solitons}
\end{equation}

\noindent Waves of the form \reff{solitons} have amplitude $3c$ and travel at constant speed $c$. 
This is contrast to solutions of linear wave equations $u_t=-cu_x$, which all travel at velocity $c$, regardless of their amplitudes.
In Fourier space, \reff{KdV} becomes
\begin{equation}
\hat{u}' = -\Dbf_N^{(3)}\hat{u} - \frac{\Dbf_N}{2}\mathbf{F}\big(\big(\mathbf{F}^{-1}\hat{u}\big)^2\big).
\end{equation}

\noindent We take
\begin{equation}
u(0,x) = 3A^2\,\mrm{sech}^2\Big(\frac{A}{2}(x+2)\Big) + 3B^2\,\mrm{sech}^2\Big(\frac{B}{2}(x+1)\Big), \quad x\in[-\pi,\pi],
\end{equation}

\noindent with $A=25$ and $B=16$, and solve up to $t=10^{-2}$. 
This is a superposition of two solitons with speed $A^2$ and $B^2$ initially centred at $x=-2$ and $x=-1$, respectively.
The stronger wave ($A=25$) catches up with the weaker one ($B=16$) at around $t=10^{-3}$.
Both waves remain unchanged after the interaction, the only nonlinear effect being a forward shift 
\begin{equation}
\frac{1}{A^2}\log\Big(\frac{A^2+B^2}{A^2-B^2}\Big)^2
\end{equation}

\noindent for the stronger wave and a backward shift 
\begin{equation}
-\frac{1}{B^2}\log\Big(\frac{A^2+B^2}{A^2-B^2}\Big)^2
\end{equation}

\noindent for the weaker one. The interaction ends at around $t=3.5\times10^{-3}$.
Figure~\ref{iconic} shows the initial condition.

\paragraph{Kuramoto--Sivashinsky} The \textit{Kuramoto--Sivashinsky equation},
\begin{equation}
u_t = - u_{xx} - u_{xxxx} - uu_x,
\label{KS}
\end{equation}

\noindent dates to the mid-1970s with the work of Kuramoto \cite{kuramoto1976} and Sivashinsky \cite{sivashinsky1977}. 
It couples destabilizing $-u_{xx}$ and stabilizing $-u_{xxxx}$ diffusions with nonlinear convection $-uu_x$. 
The nonlinear term shifts energy created at low wavenumbers by the second-order term to high wavenumbers where the fourth-order term stabilises.
The Kuramoto--Sivashinsky equation models various physical phenomena, 
from unstable drift waves in plasmas to thermal instabilities in laminar flame fronts.
In the latter, the function $u$ represents the perturbation of the flame front surface. The solutions of \reff{KS} can demonstrate a wide range of spatio-temporal dynamics, including chaos.
In Fourier space, \reff{KS} becomes
\begin{equation}
\hat{u}' = (-\Dbf_N^{(2)} - \Dbf_N^{(4)})\hat{u} - \frac{\Dbf_N}{2}\mathbf{F}\big(\big(\mathbf{F}^{-1}\hat{u}\big)^2\big).
\end{equation}

\noindent We take 
\begin{equation}
u(0,x) = \cos\Big(\frac{x}{16}\Big)\Big(1+\sin\Big(\frac{x}{16}\Big)\Big), \quad x\in[0,32\pi],
\end{equation}

\noindent and solve up to $t=100$. 
This simple initial data progressively evolves into a much more complicated superposition of wavenumbers and, even though the solution looks quite complicated, 
a characteristic pattern emerges from $t\approx50$.

\paragraph{Nonlinear Schr\"odinger} The \textit{{\normalfont (}focusing{\normalfont )} NLS equation},
\begin{equation}
u_t = iu_{xx} + i\vert u\vert^2u,
\label{NLS}
\end{equation}

\noindent models several physical phenomena, including the nonlinear propagation of light in optical fibres.
A nonlinear variant of the Schr\"odinger equation, it couples dispersion $iu_{xx}$
with a nonlinear potential $i\vert u\vert^2u$. 
Note that the wave function $u$ is complex-valued.
Among the solutions of \reff{NLS} are \textit{breathers}, given by
\begin{equation}
u(t,x) = A\bigg(\frac{2B^2\cosh(\theta)+2iB\sqrt{2-B^2}\sinh(\theta)}{2\cosh(\theta)-\sqrt{2}\sqrt{2-B^2}\cos(ABx)}-1\bigg)e^{iA^2t}, 
\end{equation}

\noindent with $\theta=A^2B\sqrt{2-B^2}t$ and $B\leq\sqrt{2}$. 
These are nonlinear waves in which energy concentrates in a localized and oscillatory fashion. 
In Fourier space, \reff{NLS} becomes
\begin{equation}
\hat{u}' = i\Dbf_N^{(2)}\hat{u} + i\mathbf{F}\Big(\big\vert\mathbf{F}^{-1}\hat{u}\big\vert^2\mathbf{F}^{-1}\hat{u}\Big).
\end{equation}

\noindent We take 
\begin{equation}
u(0,x) = \frac{2AB^2}{2-\sqrt{2}\sqrt{2-B^2}\cos(ABx)}-A, \quad x\in[-\pi,\pi],
\end{equation}

\noindent with $A=2$ and $B=1$, and solve up to $t=2$. 
This is a breather whose amplitude oscillates in time around $A=2$.
Figure~\ref{iconic} shows the initial condition.

\subsection{Model problems in 2D and 3D}

In 2D, we look for solutions of the form
\begin{equation} 
u(t,x,y) \approx \sumevenk\sumevenl \hat{u}_{k,l}(t) e^{i(kx+ly)}, \quad (x,y)\in[0,2\pi]^2,
\label{FourierSeries2D}
\end{equation}

\noindent with $N_x$ and $N_y$ points in the $x$- and $y$-directions, and appropriate rescaling for different domains. 
The unknowns, at each time $t$, are the $N_xN_y$ Fourier coefficients $\hat{u}_{k,l}(t)$.
In 3D, we look for solutions of the form
\begin{equation} 
u(t,x,y,z) \approx \sumevenk\sumevenl\sumevenm \hat{u}_{k,l,m}(t) e^{i(kx+ly+mz)}, \; (x,y,z)\in[0,2\pi]^3,
\label{FourierSeries3D}
\end{equation}

\noindent with $N_xN_yN_z$ Fourier coefficients $\hat{u}_{k,l,m}(t)$ at each time $t$.
As in 1D, the primes on the summation signs in \reff{FourierSeries2D}--\reff{FourierSeries3D} signify that the extreme terms are halved.

To construct differentiation matrices in 2D and 3D, we use Kronecker products and the 1D Fourier differentiation matrices.
For example, the Laplacian operator in 2D,
\begin{equation}
\mathcal{L}u=\Delta u=u_{xx}+u_{yy}, 
\end{equation}

\noindent is discretized by the $N_xN_y\times N_xN_y$ matrix
\begin{equation}
\Lbf=\mathbf{I}_{N_y}\otimes\mathbf{D}_{N_x}^{(2)} \, + \, \mathbf{D}_{N_y}^{(2)}\otimes\mathbf{I}_{N_x},
\label{lap2D}
\end{equation}

\noindent where $\mathbf{I}_N$ denotes the identity matrix of size $N$.

\paragraph{Ginzburg--Landau} The \textit{{\normalfont (}complex{\normalfont )} Ginzburg--Landau equation},
\begin{equation}
u_t = (1+iA)\Delta u + u - (1+ iB)u\vert u\vert^2,
\label{GL}
\end{equation}

\noindent was first derived in 2D by Stewartson and Stuart in 1971 to study nonlinear instabilities in plane Poiseuille flow \cite{stewartson1971}, using concepts from Ginzburg--Landau theory for superconductivity.
The function $u$ is the amplitude of a nonlinear perturbation wave for values of the Reynolds number close to the critical value, above which perturbations may grow.
Equation \reff{GL} admits solutions known as \textit{frozen states} which correspond to quasi-frozen spiral defects surrounded by shock lines. In this regime, $\vert u\vert$ is stationary in time.
We take $A=0$ and $B=1.5$,
\begin{equation}
u(0,x,y) = e^{-0.1[(x-50)^2 + (y-50)^2]}, \quad (x,y)\in[0,100]^2,
\end{equation}

\noindent and 
\begin{equation}
u(0,x,y,z) = e^{-0.1[(x-50)^2 + (y-50)^2 + (z-50)^2]},  \quad (x,y,z)\in[0,100]^3,
\end{equation}

\noindent and solve up to $t=10$ in both 2D and 3D. These two initial conditions generate spiral waves.
Figure~\ref{iconic} shows an example of a 2D frozen state solution at $t=30$, obtained with an initial condition of amplitude $0.1$ involving random noise on the grid.

\paragraph{Schnakenberg} The \textit{Schnakenberg equations},
\begin{equation}
\left\{
\begin{array}{l}
u_t = \epsilon_u\Delta u + \gamma(a - u + u^2v), \\\\
v_t = \epsilon_v\Delta v + \gamma(b - u^2v),
\end{array}
\right. 
\label{Schnak}
\end{equation}

\noindent are reaction-diffusion equations derived by Schnakenberg in 1979 to study limit cycle behaviours of two-component chemical reactions \cite{schnakenberg1979}.
The system \reff{Schnak} models the chemical reaction $2U+V\rightarrow3U\,;\,U\rightleftharpoons A\,;\,B\rightarrow V$; $u$ and $v$ are the concentrations of $U$ and $V$,
and $A$ and $B$ are another two chemical species, assumed to be maintained at constant concentrations $a$ and $b$.
We take $\epsilon_u=1$, $\epsilon_v=10$, $\gamma=3$, $a=0.1$ and $b=0.9$. The initial conditions are
\begin{equation}
\begin{array}{l}
u(0,x,y) = 1 - e^{-2[(x-G/2.15)^2 + (y-G/2.15)^2]}, \\\\
\dsp v(0,x,y) = \frac{0.9}{(0.1+0.9)^2} + e^{-2[(x-G/2)^2 + 2(y-G/2)^2]},
\end{array}
\end{equation}

\noindent with $(x,y)\in[0,G]^2$ and $G=30$ in 2D, and 
\begin{equation}
\begin{array}{l}
u(0,x,y,z) = 1 - e^{-2[(x-G/2.15)^2 + (y-G/2.15)^2 + (z-G/2.15)^2]}, \\\\
\dsp v(0,x,y,z) = \frac{0.9}{(0.1+0.9)^2} + e^{-2[(x-G/2)^2 + 2(y-G/2)^2 + 2(z-G/2)^2]},
\end{array}
\end{equation}

\noindent with $(x,y,z)\in[0,G]^3$ and $G=30$ in 3D. 
We solve up to $t=20$ in both 2D and 3D.
Note that these initial conditions are small perturbations from the constant solution $(u,v)=(a+b,b/(a+b)^2)$.
They lead to a set of spots at around $t=500$ in 2D and $t=300$ in 3D.
Figure~\ref{iconic} shows the 2D solution at $t=500$.

\paragraph{Swift--Hohenberg} The \textit{Swift--Hohenberg equation},
\begin{equation}
u_t = ru - (1+\Delta)^2u + gu^2-u^3,
\label{SH}
\end{equation}

\noindent was first derived in 2D by Swift and Hohenberg in 1977 to study thermal fluctuations on a fluid near the Rayleigh--B\'{e}nard convective instability \cite{swift1977}. In 2D, the function $u$ is the temperature field in a plane horizontal layer of fluid heated from below.
Equation \reff{SH} is another example of a PDE that exhibits pattern formation, including stripes, spots and spirals.
We take $r=0.1$, $g=1$, 
\begin{equation}
u(0,x,y) = \frac{1}{4}\Big(\sin\Big(\frac{\pi x}{10}\Big) + \sin\Big(\frac{\pi y}{10}\Big) + \sin\Big(\frac{\pi x}{2}\Big)\sin\Big(\frac{\pi y}{2}\Big)\Big),
\end{equation}
 
\noindent with $(x,y)\in[0,20]^2$ and 
\begin{equation}
\begin{array}{l}
u(0,x,y,z) = \dsp\frac{1}{4}\Big(\sin\Big(\frac{\pi x}{10}\Big) + \sin\Big(\frac{\pi y}{10}\Big) + \sin\Big(\frac{\pi z}{10}\Big) + \sin\Big(\frac{\pi x}{2}\Big)\sin\Big(\frac{\pi y}{2}\Big)\\\\
\dsp\hspace{2.5cm} + \sin\Big(\frac{\pi x}{2}\Big)\sin\Big(\frac{\pi z}{2}\Big) + \sin\Big(\frac{\pi y}{2}\Big)\sin\Big(\frac{\pi z}{2}\Big)\Big),
\end{array}
\end{equation}

\noindent with $(x,y,z)\in[0,20]^3$, and solve up to $t=20$ in both 2D and 3D. Both of these two initial conditions lead to a set of spots.
Figure~\ref{iconic} shows the 2D solution at $t=1000$ obtained with $r=0.1$, $g=0$ and an initial condition of amplitude $0.1$ involving random noise on the grid.
This solution corresponds to the so-called \textit{convection rolls}.

\section{Numerical comparisons}

\subsection{Methodology}

\begin{figure}
\centering
\includegraphics [scale=.35]{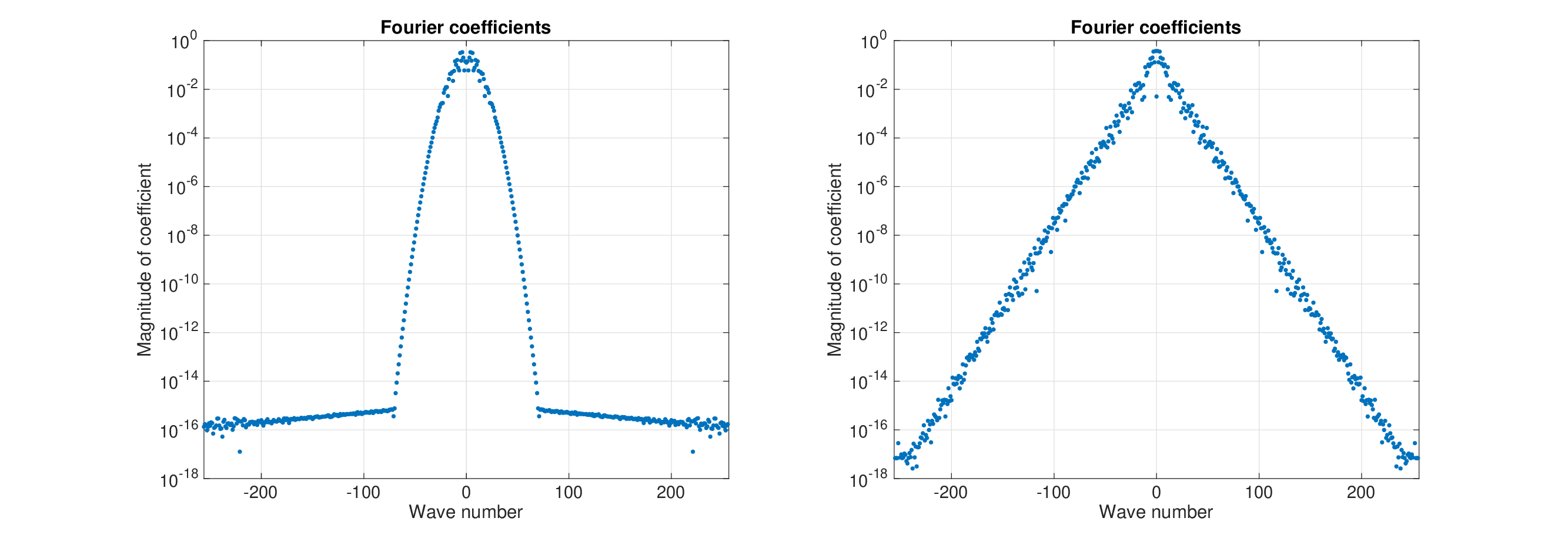}
\caption{\textit{Fourier coefficients of the initial condition $\reff{ACIC}$ of the Allen--Cahn equation (left) and of the solution at $t=60$ computed
with time-step $h=10^{-2}$ (right). With $N=512$ grid points, the Fourier coefficients decay to about $10^{-16}$. In our computations for the Allen--Cahn equation, 
the smallest error due to the time discretization is equal to about $10^{-12}\gg10^{-16}$.}}
\label{Fourier}
\end{figure}

To compare exponential integrators, we follow the methodology of~\cite{kassam2005}, though the experiments described here are far more extensive.
We solve a given PDE up to $t=T$ for various time-steps $h$ and a fixed number of grid points.
We estimate the ``exact'' solution $u^{ex}(t=T,X)$ by using a ``very small'' time-step (half the smallest time-step $h$)
and the PECEC736 scheme (one of the two seventh-order accurate schemes in Table~\ref{ReferenceTable}).
We then measure the relative $L^2$-error $E$ at $t=T$ between the computed solution $u(t=T,X)$ and $u^{ex}(t=T,X)$, i.e.,
\begin{equation}
E =\frac{\Vert u(t=T,X) - u^{ex}(t=T,X)\Vert_2}{\Vert u^{ex}(t=T,X)\Vert_2}.
\label{error}
\end{equation}

\noindent For both $u$ and $u^{ex}$ we use $N=512$ grid points in 1D, $N_x=N_y=128$ grid points in 2D and $N_x=N_y=N_z=128$ grid points in 3D.
(With these grid sizes, the error due to the spatial discretization is small compared to the error due to the time discretization; see Figure~\ref{Fourier}.)
For the contour integrals, we use $M=64$ points in 1D and $M=32$ points in 2D and 3D.
We plot \reff{error} against relative time-steps $h/T$ and computer times on a pair of graphs.\footnote{The precomputation of the coefficients of the exponential integrators and the starting phase of multistep methods are not included in the computing time. 
Timings were done on a 2.8\,GHz Intel i7 machine with 16\,GB of RAM.}
The former gives a measure of the accuracy of the exponential integrator for various time-steps or, equivalently, for various number of integration steps.
(If the relative time-step is $10^{-3}$, it means that the integrator performed $10^3$ steps to reach $t=T$.) 
However, it is possible that each step is more costly, so it is the latter that ultimately matters. 
We compare different families on different pairs of graphs with curves for ETDRK4 included on all plots as a baseline. 
We have tested every integrator on every PDE, but we shall only show graphs that correspond to the characteristic behaviour of a family of integrators,
or highlight notable features such as instability or particularly good/bad performance.
The rest of the graphs can be found in the Ph.D. thesis of the first author \cite{montanelli2017}.

\subsection{Starting multistep schemes}

To start a multistep scheme with $q$ steps, one needs $q$ values: the initial condition $\hat{u}^0$ and $q-1$ extra values 
$\hat{u}^1,\ldots,\hat{u}^{q-1}$.
It is suggested in~\cite{calvo2006} to use the following strategy: first, compute an approximation of $\hat{U}=(\hat{u}^1,\ldots,\hat{u}^{q-1})^T$ using 
a low-order exponential integrator (e.g., ETDRK2, the second-order version of ETDRK4, also introduced by Cox and Matthews in~\cite{cox2002}), and then,
use a fixed point iteration to refine this approximation. The fixed point iteration is applied to the following system of nonlinear equations,
\begin{equation}
\hat{u}^j = e^{jh\Lbf}\hat{u}^0 + h\sum_{l=0}^{q-1} \gamma_l(j, h\Lbf)\Delta^l \Nbf(\hat{u}^0), \quad 1\leq j\leq q-1,
\label{nonlinsys}
\end{equation}

\noindent where $\Delta^l$ is the forward difference operator,
\begin{equation}
\Delta^0 \Nbf(\hat{u}^0) = \hat{u}^0, \quad \Delta^l \Nbf(\hat{u}^0) = \sum_{i=0}^l (-1)^i {l\choose i} \Nbf(\hat{u}^{l-i}), \quad l\geq 1,
\end{equation}

\noindent and the $\gamma$-functions are defined by the recurrence relation,
\begin{equation}
\begin{array}{l}
\dsp \gamma_0(k,z)=\frac{e^{kz}-1}{z}, \\\\
\gamma_j(k, z) = \frac{\dsp \Bigg(\sum_{m=1}^{j} \frac{(-1)^{m-1}}{m}\gamma_{j-m}(k,z)\Bigg) - {k\choose j}}{\dsp z}, \quad 1\leq j\leq k, \\\\
\gamma_j(k, z) = \frac{\dsp \Bigg(\sum_{m=1}^{j} \frac{(-1)^{m-1}}{m}\gamma_{j-m}(k,z)\Bigg)}{\dsp z}, \quad j > k.
\end{array}
\end{equation}

\noindent Note that, like the $\varphi$-functions, the $\gamma$-functions can be evaluated by contour integrals and satisfy
the symmetry property $\gamma(\bar{z})=\bar{\gamma}(z)$ for all $z\in\C$.
Let us write~\reff{nonlinsys} as $\hat{U} = F(\hat{U})$. The fixed point iteration is then given by
\begin{equation}
\hat{U}_{[n+1]} = F(\hat{U}_{[n]}),
\label{iteration}
\end{equation}

\noindent where $\hat{U}_{[n]}=(\hat{u}^1_{[n]},\ldots,\hat{u}^{q-1}_{[n]})^T$ denotes the approximation obtained
after $n$ iterations ($\hat{U}_{[0]}$ corresponding to the approximation given by ETDRK2).
The fixed point iteration~\reff{iteration} is carried out until the norm of the difference between two successive iterations is of the order of $h^q$.

\subsection{Results}

\begin{figure}
\hspace{-1.2cm}
\includegraphics [scale=.4]{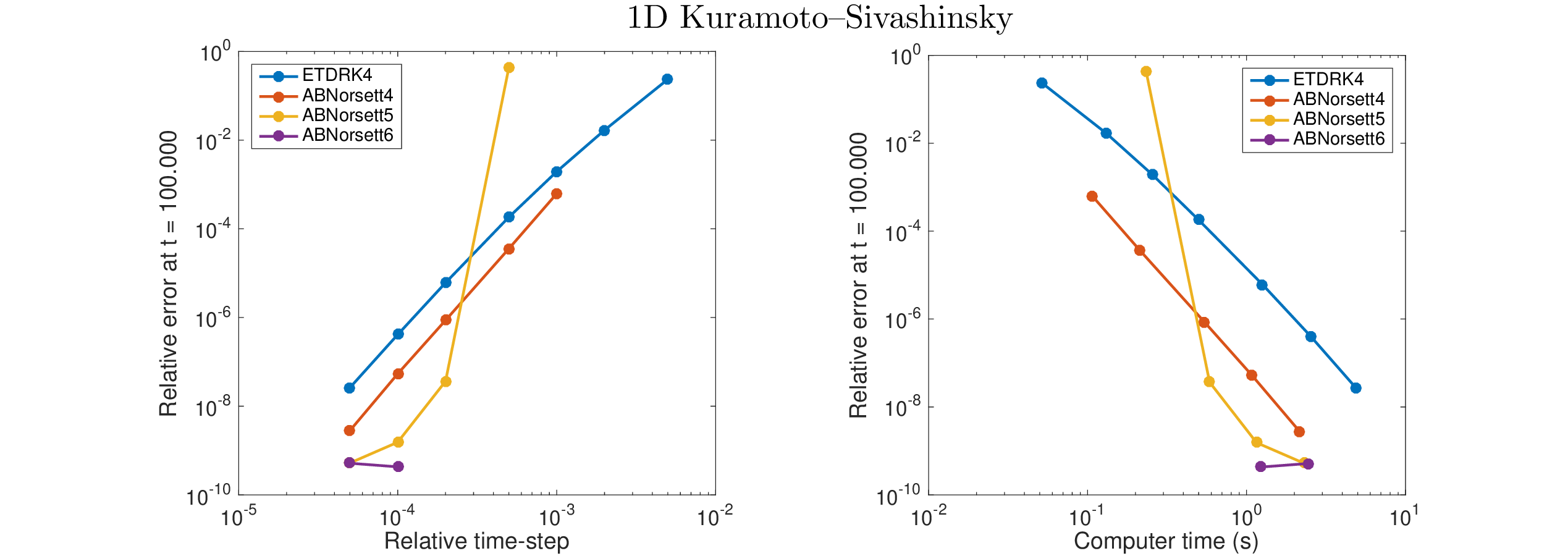}
\caption{\textit{Accuracy versus time-step and computer time for the ETD Adams--Bashforth methods for the $1$D Kuramoto--Sivashinsky equation.
The order $5$ and $6$ methods are more efficient than ETDRK$\,4$ at high accuracies, but often unstable at lower accuracies, as reflected in dots missing from the curves.}}
\label{ETDAB_graphs}
\end{figure}

\begin{figure}
\hspace{-1.2cm}
\includegraphics [scale=.4]{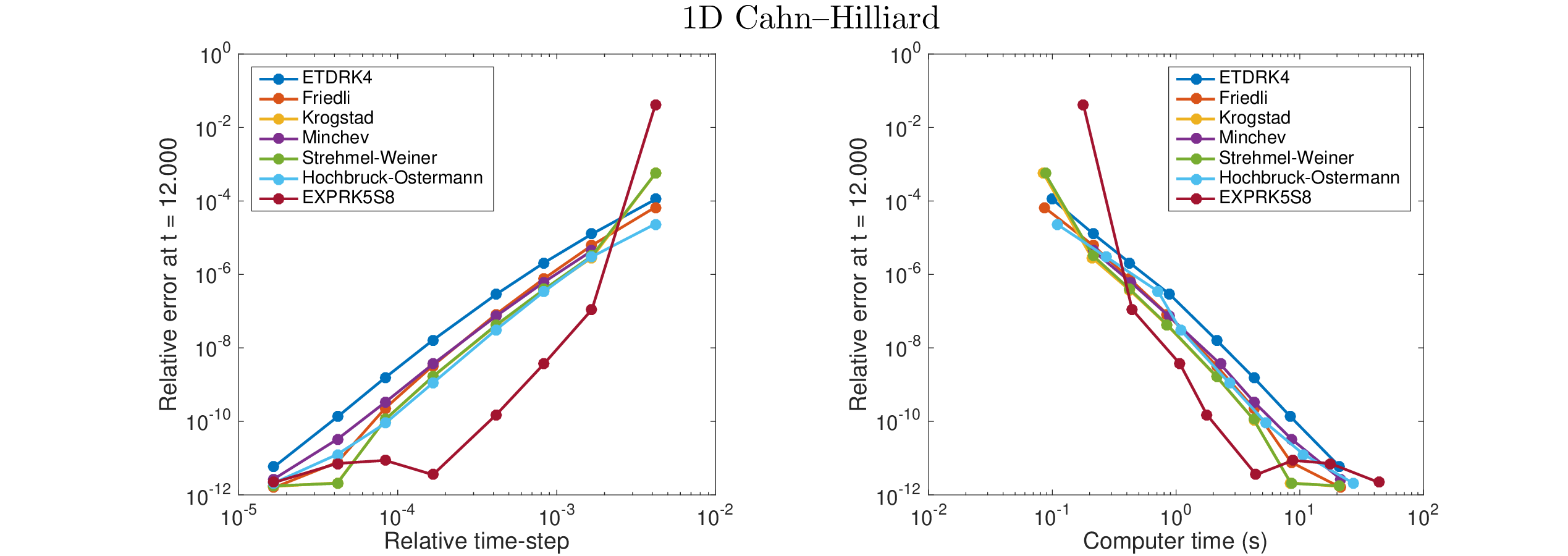}
\par\hspace{-1.2cm}
\includegraphics [scale=.4]{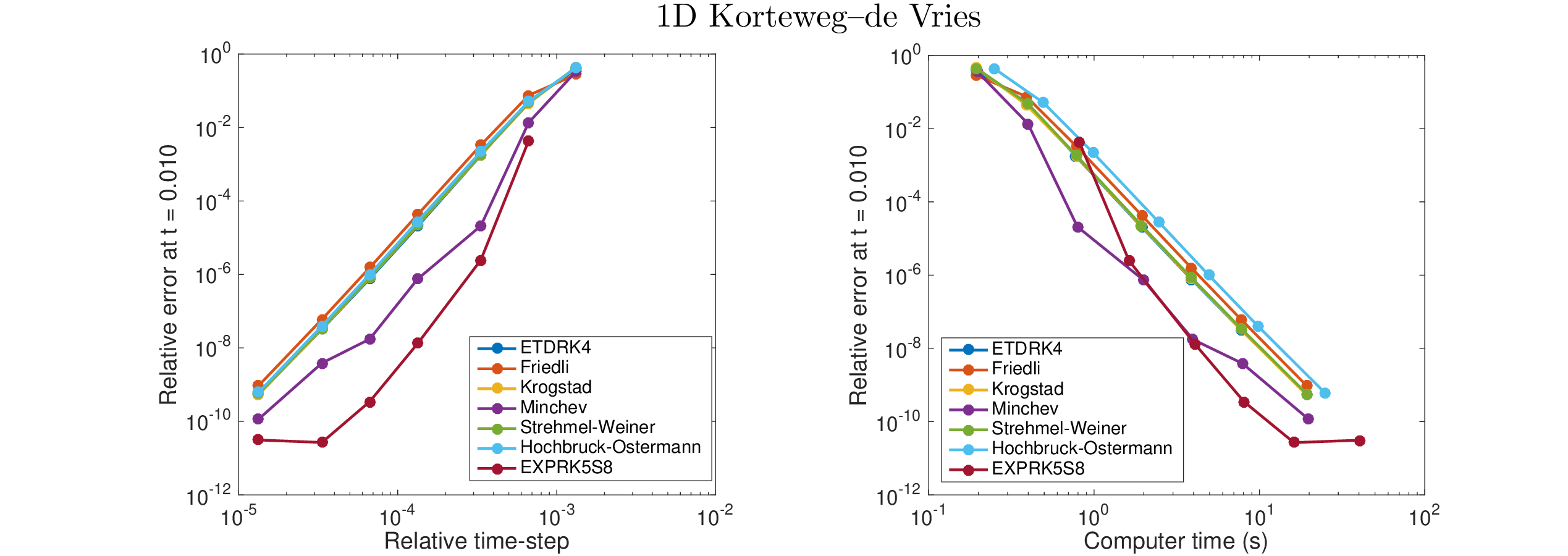}
\par\hspace{-1.2cm}
\includegraphics [scale=.4]{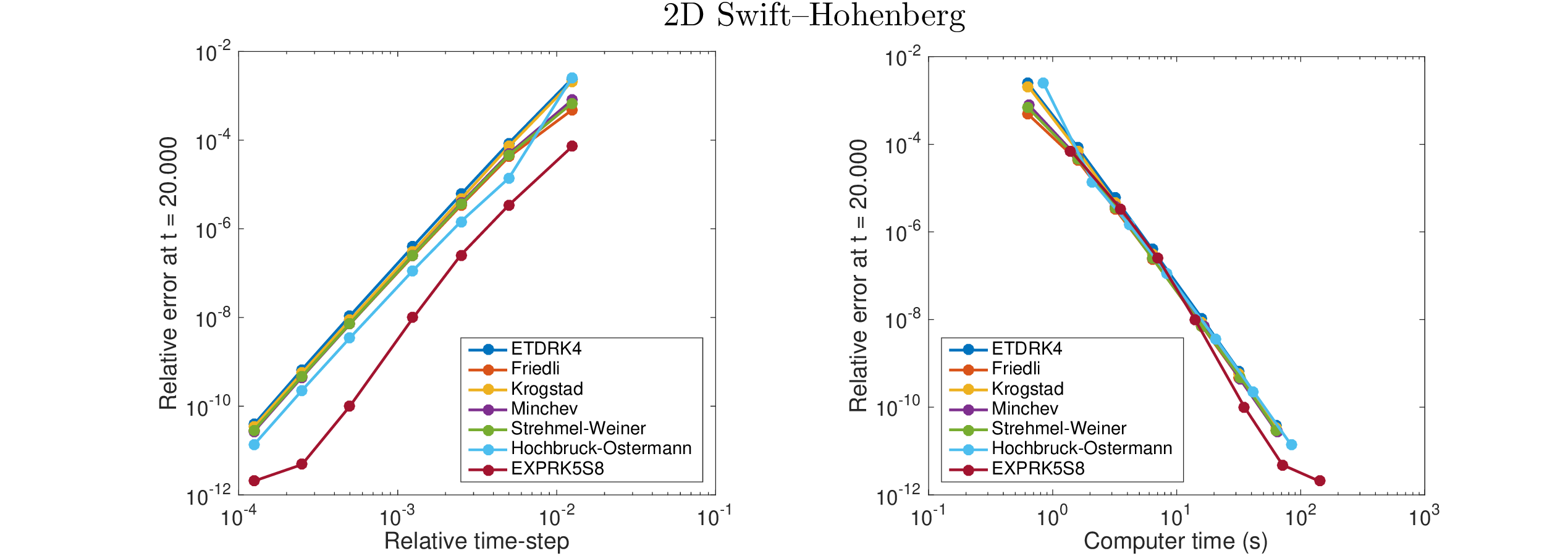}
\caption{\textit{Accuracy versus time-step and computer time for the ETD Runge--Kutta methods for the $1$D Cahn--Hilliard (top), $1$D KdV (centre) and $2$D Swift--Hohenberg (bottom) equations.
The EXPRK$\,5$S$\,8$ scheme is impressively efficient in $1$D but is unstable at low accuracies for the Cahn--Hilliard and KdV equations.
In $2$D and $3$D, it does not beat the fourth-order methods for the Schnakenberg (Figures {\nf C.13} and {\nf C.19} in {\nf \cite{montanelli2017}}) and Swift--Hohenberg equations (above for $2$D, Figure {\nf C.21} in {\nf \cite{montanelli2017}} for $3$D).}}
\label{ETDRK_graphs}
\end{figure}

\begin{figure}
\hspace{-.8cm}
\includegraphics [scale=.37]{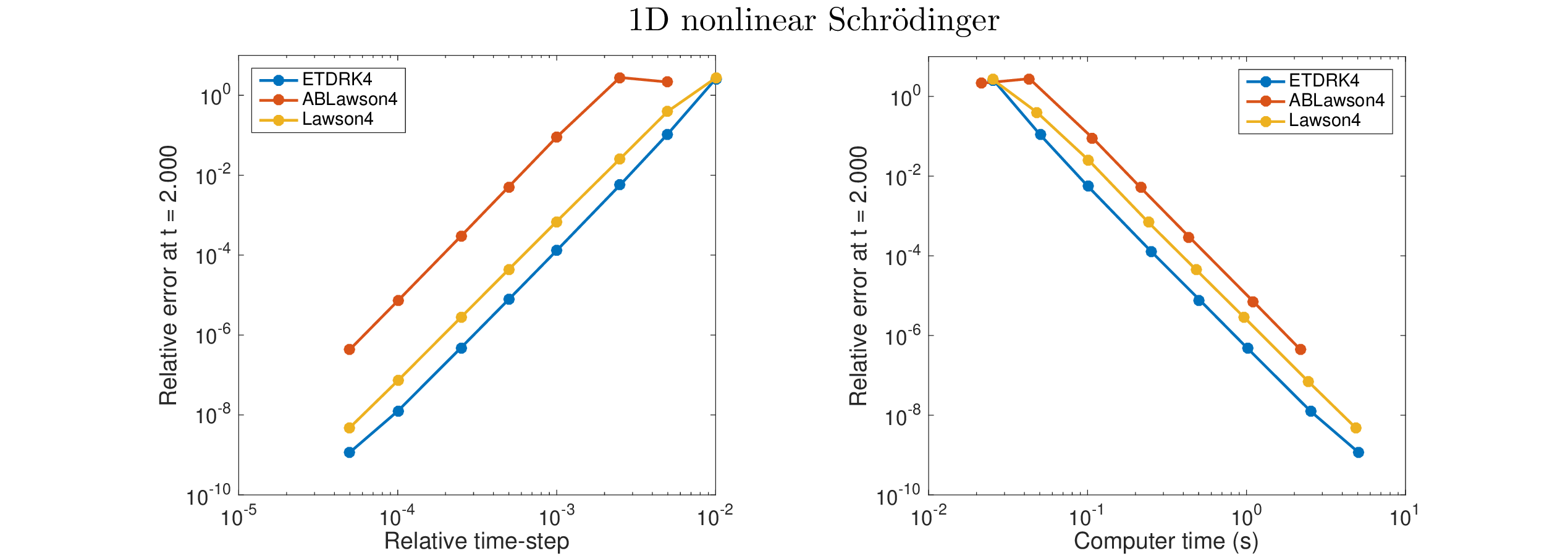}
\caption{\textit{Accuracy versus time-step and computer time for the Lawson methods for the $1$D NLS equation.
These formulas are too inaccurate to be competitive; the constants involved in the convergence bounds are too great.}}
\label{Lawson_graphs}
\end{figure}

\begin{figure}
\hspace{-.8cm}
\includegraphics [scale=.37]{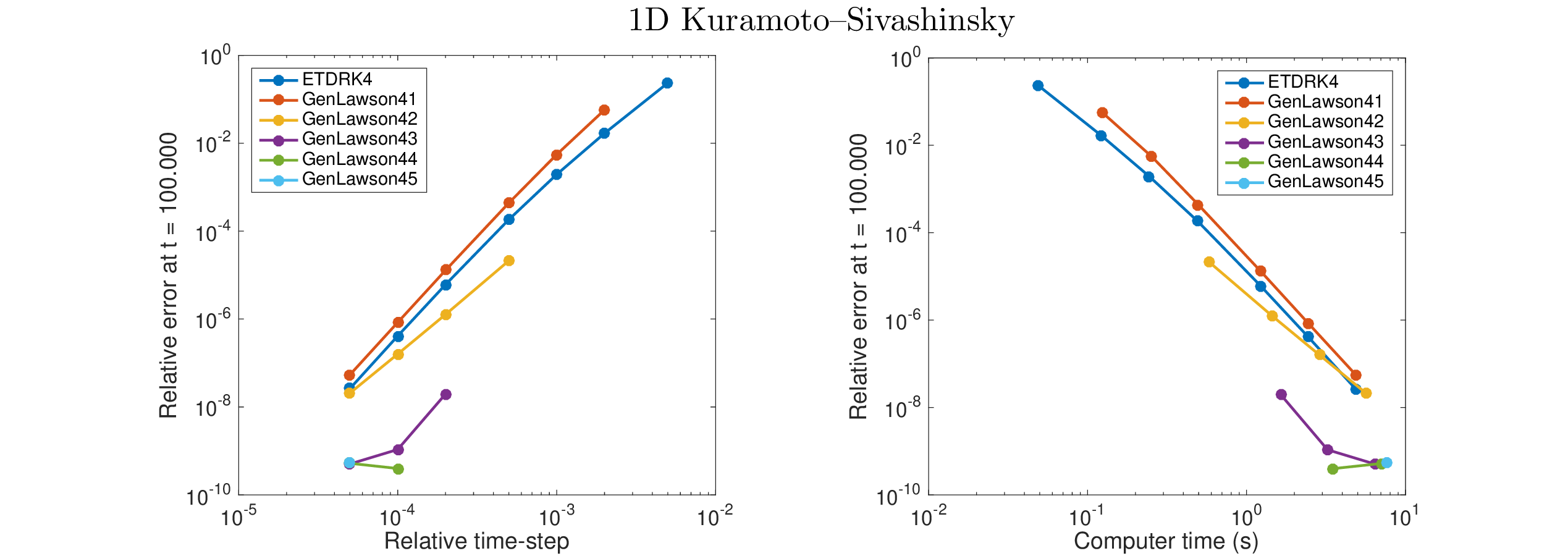}
\par\hspace{-.8cm}
\includegraphics [scale=.37]{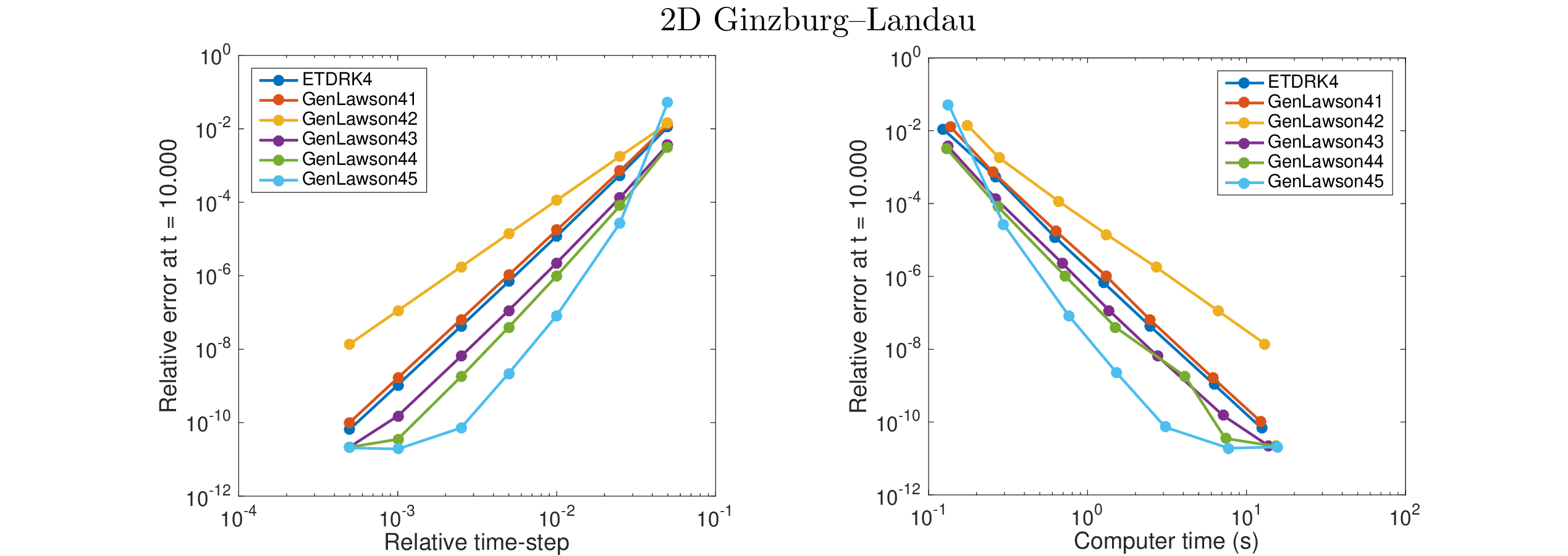}
\caption{\textit{Accuracy versus time-step and computer time for the generalised Lawson methods for the $1$D Kuramoto--Sivashinsky (top) and $2$D Ginzburg--Landau (bottom) equations.
These methods are highly unstable for most PDEs in $1$D. 
In $2$D and $3$D, the GenLawson$\,41$ formula has virtually identical performance to ETDRK$\,4$ while the other variants with three to five steps perform well for some problems (e.g., $2$D Ginzburg--Landau equation above) but are less efficient for others (e.g., Figure {\nf C.14} in {\nf \cite{montanelli2017}}).}}
\label{GenLawson_graphs}
\end{figure}

\begin{figure}
\hspace{-1.25cm}
\includegraphics [scale=.4]{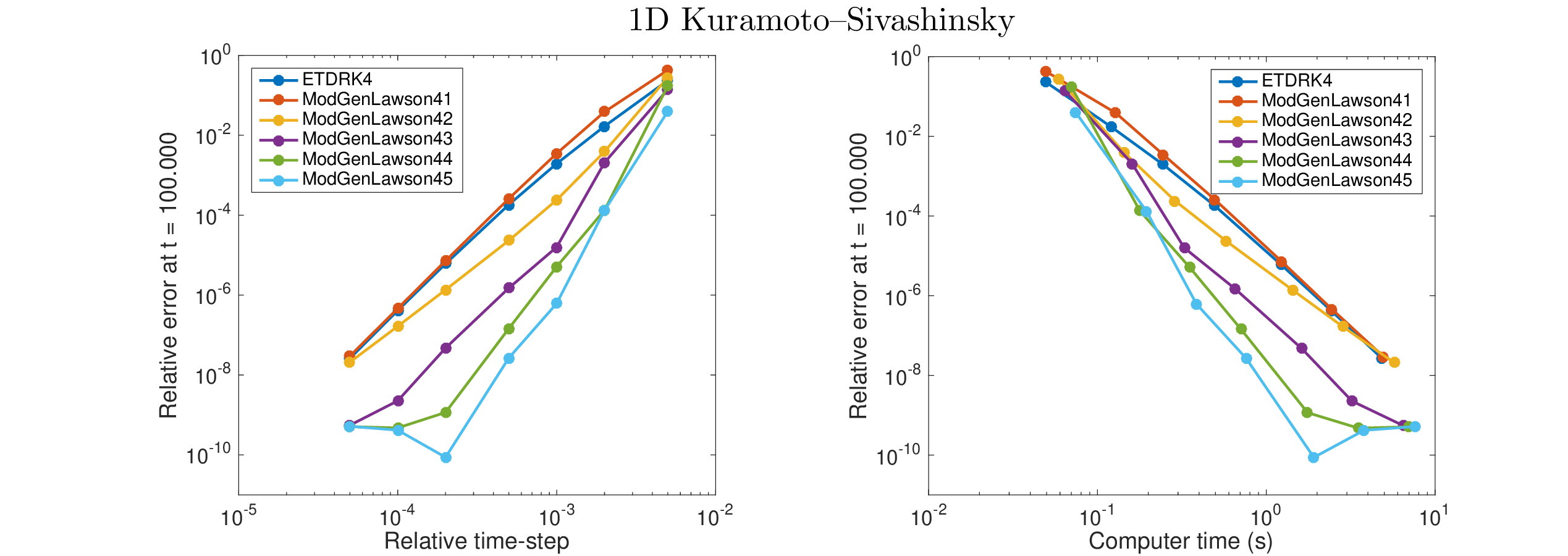}
\par\hspace{-1.25cm}
\includegraphics [scale=.4]{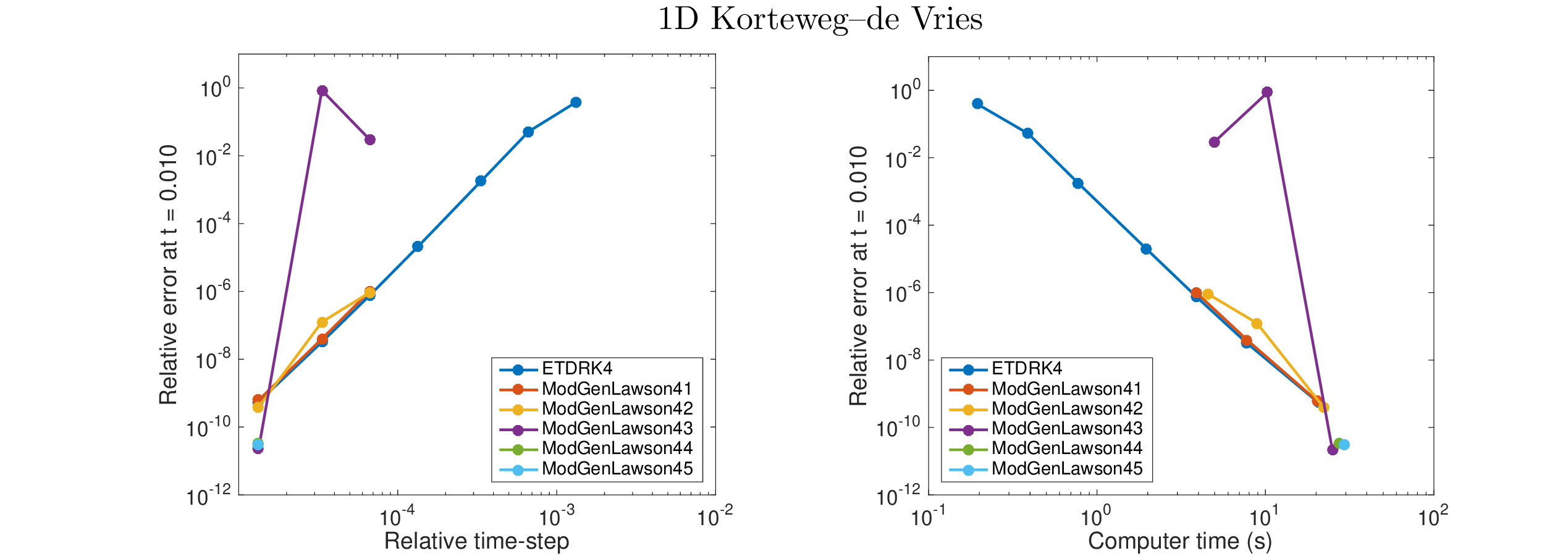}
\par\hspace{-1.25cm}
\includegraphics [scale=.4]{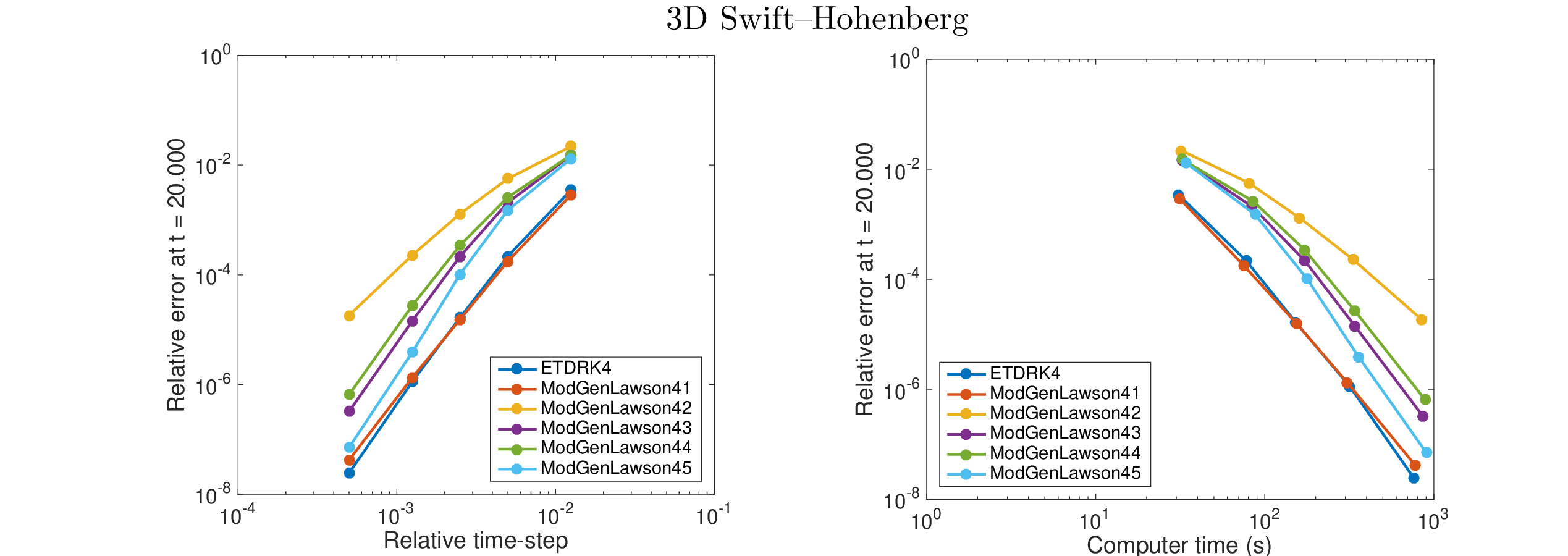}
\caption{\textit{Accuracy versus time-step and computer time for the modified generalised Lawson methods for the $1$D Kuramoto--Sivashinsky (top), $1$D KdV (centre) and $3$D Swift--Hohenberg (bottom) equations.
In $1$D, these methods are much more stable than the generalised Lawson methods but are still highly unstable for the KdV equation.
In $2$D and $3$D, they are very similar to the generalised Lawson methods, i.e., they perform well for some problems  (e.g., Figure {\nf C.18} in {\nf \cite{montanelli2017}}) but are less efficient for others  (e.g., $3$D Swift--Hohenberg equation above).}}
\label{ModGenLawson_graphs}
\end{figure}

\begin{figure}
\hspace{-1.1cm}
\includegraphics [scale=.4]{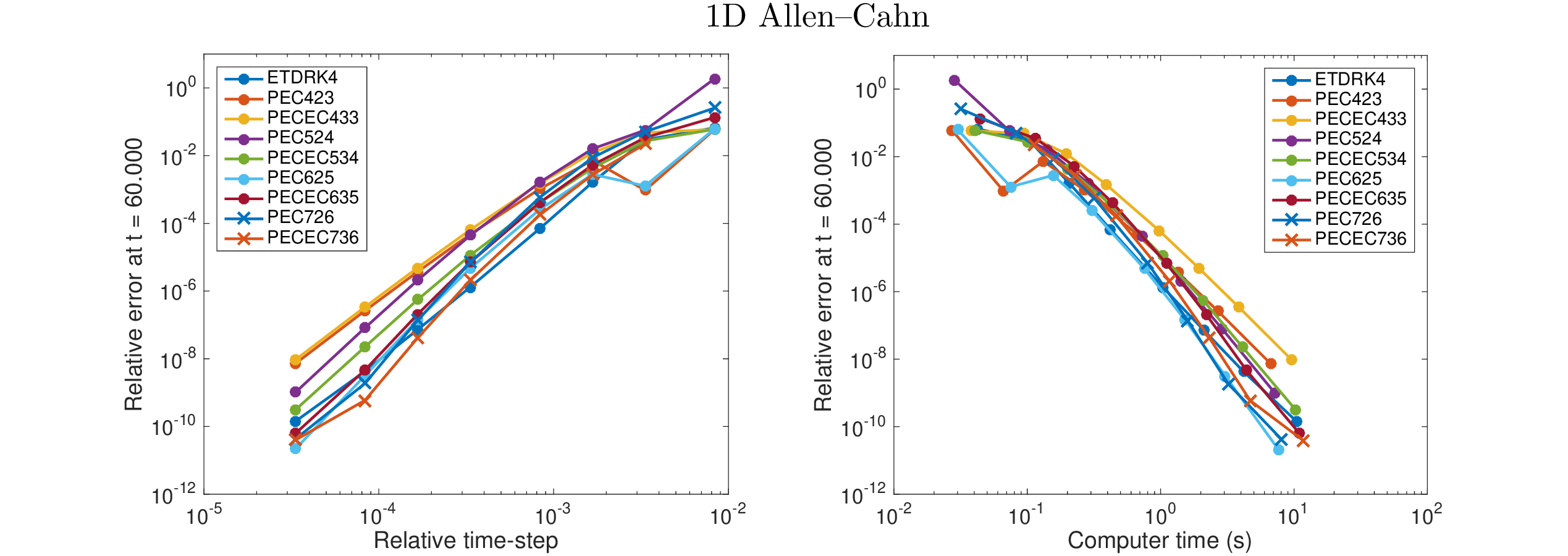}
\par\hspace{-1.1cm}
\includegraphics [scale=.4]{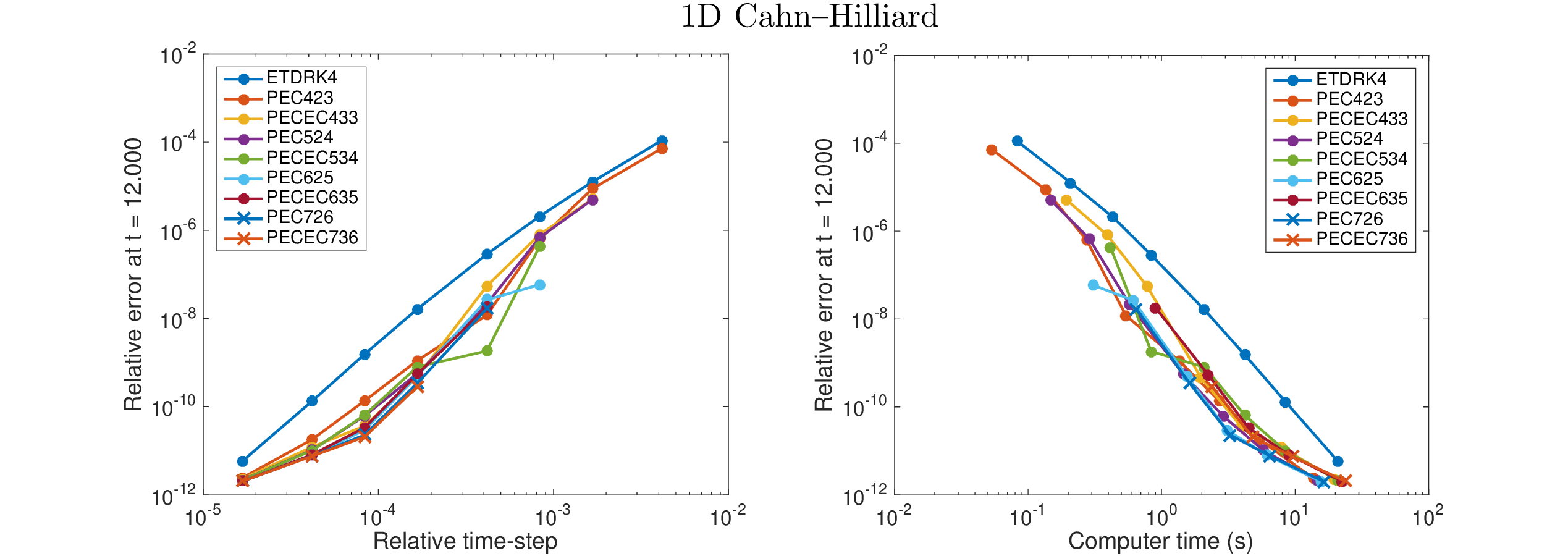}
\par\hspace{-1.1cm}
\includegraphics [scale=.4]{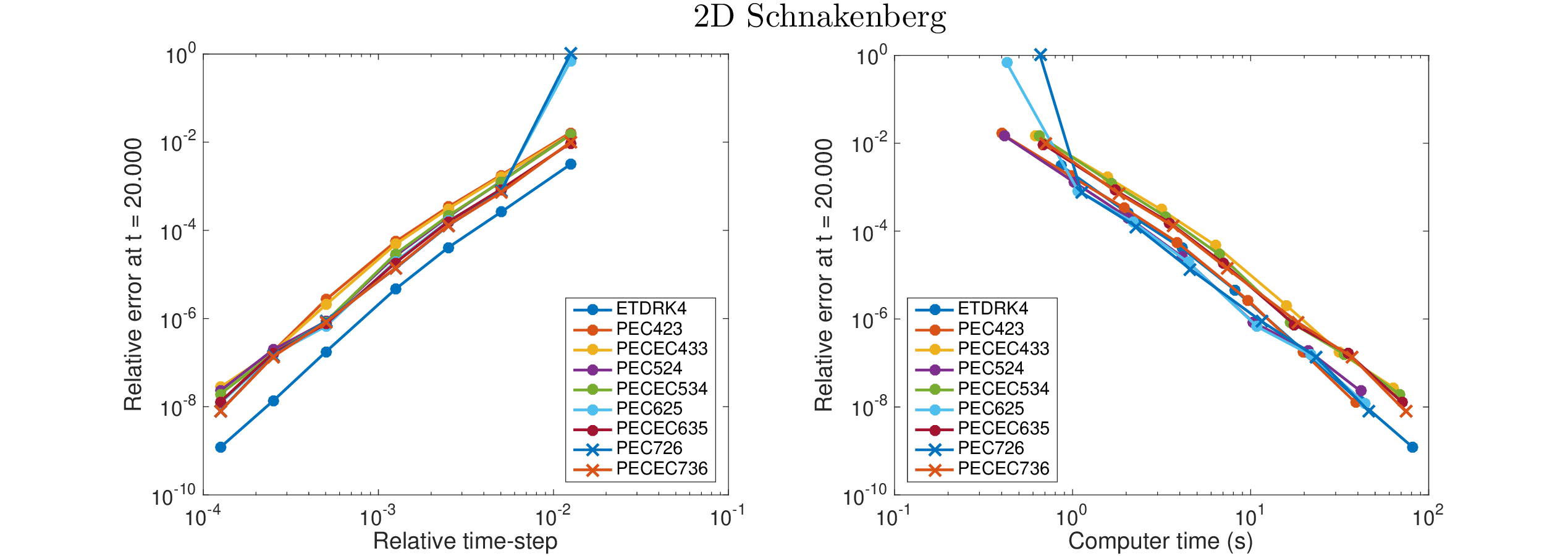}
\caption{\textit{Accuracy versus time-step and computer time for the predictor-corrector methods for the $1$D Allen--Cahn (top),
$1$D Cahn--Hilliard (centre) and $2$D Schnakenberg equations (bottom).
In $1$D, these schemes are efficient, but do not clearly outperform ETDRK$\,4$ for the Allen--Cahn equation and are unstable at low accuracies for the Cahn--Hilliard equation.
In $2$D and $3$D, they beat ETDRK$\,4$ for the Ginzburg--Landau equation (Figures {\nf C.12} and {\nf C.18} in {\nf \cite{montanelli2017}}) but have similar performance for the other equations (e.g., $2$D Schnakenberg equations above).}}
\label{PEC_graphs}
\end{figure}

We now report on the results of our numerical testing and present plots showing the typical behaviours and notable features we see from each set of methods.

Figure~\ref{ETDAB_graphs} shows results for the \textit{ETD Adams--Bashforth} methods for the Kuramoto--Sivashinsky equation (1D). 
These formulas are often unstable for large time-steps but can be competitive at high accuracies.
For the Cahn--Hilliard and KdV equations (Figures C.3 and C.5 in \cite{montanelli2017}), we could not get them to work at all with the spatial discretization that we used.

Figure~\ref{ETDRK_graphs} shows results for the \textit{ETD Runge--Kutta} methods for the 1D Cahn--Hilliard, 1D KdV and 2D Swift--Hohenberg equations.
These formulas have good stability properties.
The fourth-order methods have similar performance to ETDRK4.
The fifth-order EXPRK5S8 integrator is impressively efficient in 1D, but exhibits instability for the Cahn--Hilliard and KdV equations.
In 2D and 3D, it is more accurate than the fourth-order methods for the six PDEs we have considered, and more efficient for the Ginzburg--Landau equation (Figures C.11 and C.17 in \cite{montanelli2017}).
However, for the Schnakenberg (Figures C.13 and C.19 in \cite{montanelli2017}) and Swift--Hohenberg equations, it is not accurate enough to counterbalance its high cost per time-step.

Figure~\ref{Lawson_graphs} shows results for the \textit{Lawson} methods for the NLS equation (1D).
These formulas are not accurate enough to be competitive.
For the ABLawson4 formula, this lack of accuracy is partly compensated by its low computational cost per time-step (it is a purely multistep method).

Figure~\ref{GenLawson_graphs} shows results for the \textit{generalised Lawson} methods for the 1D Kuramoto--Sivashinsky and 2D Ginzburg--Landau equations.
In 1D, these formulas suffer from instabilities for most PDEs.
In 2D and 3D, the GenLawson41 formula has virtually identical performance to ETDRK4, the GenLawson42 formula is always less efficient than ETDRK4, while the other variants with three to five steps perform well for the Ginzburg--Landau equation but are less efficient for the Schnakenberg and Swift--Hohenberg equations (Figures C.14, C.16, C.20 and C.22 in \cite{montanelli2017}).

Figure~\ref{ModGenLawson_graphs} shows results for the \textit{modified generalised Lawson} methods for the 1D Kuramoto--Sivashinsky, 1D KdV and 3D Swift--Hohenberg equations.
In 1D, these formulas are much more stable than the generalised Lawson schemes and are quite efficient, but still suffer from instabilities for the 1D KdV equation.
In 2D and 3D, we reach the same conclusions as for the generalised Lawson methods: the ModGenLawson41 formula has virtually identical performance to ETDRK4, the ModGenLawson42 formula is always the least efficient and the other variants perform well for some problems but are less efficient for others.

Figure~\ref{PEC_graphs} shows results for the \textit{exponential predictor-corrector} methods for the 1D Allen--Cahn, 1D Cahn--Hilliard and 2D Schnakenberg equations.
These formulas are particularly efficient in 1D, especially for the Kuramoto--Sivashinsky and NLS equations (Figures C.8 and C.10 in \cite{montanelli2017}), but do not clearly outperform ETDRK4 for the Allen--Cahn equation.
Most of them are unstable at low accuracies for the Cahn--Hilliard equation, especially the higher-order schemes, and some of them are also unstable at low accuracies for the KdV equation (Figure C.6 in \cite{montanelli2017}).
In 2D and 3D, they are more efficient than ETDRK4 for the Ginzburg--Landau equation (Figures C.12 and C.18 in \cite{montanelli2017}) but most of them have
similar performance to it for the other PDEs we have considered.
Note that the higher-order schemes with two steps (PEC625 and PEC726) also beat ETDRK4 for the Swift--Hohenberg equations, 
but these are particularly unstable for the Cahn--Hilliard and KdV equations.

\section{Discussion}

We have tested 30 exponential integrators on 11 model problems in 1D, 2D and 3D, and have observed considerable differences in stability and efficiency.
As expected, the schemes did not exhibit any order reduction (periodic boundary conditions).
The main conclusion is that it is difficult to find a method that outperforms ETDRK4 for all the PDEs we have considered.

Our experiments show that the ETD Adams--Bashforth and the generalised Lawson methods are highly unstable while the Lawson methods are not accurate enough.
Within  the ETD Runge--Kutta methods, it is hard to do much better than ETDRK4. 
The fourth-order schemes are quite similar in terms of efficiency and stability.
The fifth-order EXPRK5S8 integrator is more efficient than ETDRK4 for most PDEs in 1D, but is unstable at low accuracies for the KdV and Cahn--Hilliard equations.
In 2D and 3D, it outperforms ETDRK4 only for the Ginzburg--Landau equation.
Since it requires the precomputation of more than twice as many coefficients as ETDRK4, it makes it much more complicated to implement and probably less appealing to general users.
The high-order modified generalised Lawson and exponential predictor-corrector methods are competitive stiff solvers for some PDEs, but for others do not outperform ETDRK4 or else suffer from instabilities.

Our numerical experiments were performed using MATLAB and have been embedded within Chebfun.
More specifically, the \texttt{spin}, \texttt{spin2} and \texttt{spin3} codes implement a Fourier spectral method and exponential integrators to solve PDEs in 1D, 2D and 3D periodic domains.
(Note that \texttt{spin} stands for \textbf{s}tiff \textbf{P}DE \textbf{in}tegrator.)
These have been one of the most major additions to Chebfun in recent years from a user point of view.
The simplest way to see \texttt{spin} in action is to type simply \texttt{spin('ks')} (for the Kuramoto--Sivashinsky equation) or \texttt{spin2('gl2')} (for the 2D Ginzburg--Landau equation) to invoke an example computation.
It is also possible to define your own PDE using the \texttt{spinop} class.
To produce the graphs of Section 4.3, we have used the \texttt{spincomp} code.

\section*{Acknowledgements}

This paper is dedicated to Nick Trefethen for his inspirational contributions to the field of numerical analysis. 

\bibliographystyle{siam}
\bibliography{spin}

\end{document}